\documentclass[a4paper,11pt]{article}
\pdfoutput=1
\usepackage{mathrsfs}
\usepackage[utf8]{inputenc}
\usepackage[T1]{fontenc}
\usepackage{amsmath}
\usepackage{amssymb}
\usepackage{amsthm}
\usepackage[francais,english]{babel}
\usepackage{graphicx}
\usepackage{color}
\usepackage{fullpage}
\usepackage{dsfont}
\usepackage{enumerate}
\usepackage{caption}
\usepackage{bbm}

\definecolor{purple}{rgb}{0.57,0.1,0.53}

\newcommand{\HypLevy}{A}
\newcommand{\HypMut}{B.2}
\newcommand{\HypMutConst}{B.1}

\newtheorem{lemme}{Lemma}[section]
\newtheorem{theoreme}[lemme]{Theorem}

\newtheorem{proposition}[lemme]{Proposition}
\newtheorem{corollaire}[lemme]{Corollary}

\newtheorem{remarque}[lemme]{Remark}
\newenvironment{demo}{\noindent\emph{Proof :} \\ }{\hfill $\square$ \\ \par\smallskip}
\newenvironment{demoth}[1]{\noindent\emph{Proof of Theorem #1 :} \\ }{\hfill $\square$}
\newenvironment{demopr}[1]{\noindent\emph{Proof of Proposition #1 :} \\ }{\hfill $\square$}
\newenvironment{democo}[1]{\noindent\emph{Proof of Corollary #1 :} \\ }{\hfill $\square$}

\newcounter{claim}

\newtheorem*{theoreme*}{Theorem}
\newtheorem*{proposition*}{Proposition}
\newtheorem*{remarque*}{Remark}

\newcommand{\Z}{\mathbb Z}
\renewcommand{\L}{\Lambda}

\newcommand{\sachant}{\,|\,}
\newcommand{\R}{\mathbb R}
\newcommand{\N}{\mathbb N}

\renewcommand{\d}[1]{\text{d}#1}
\renewcommand{\P}{\mathbb P}
\newcommand{\E}{\mathbb E}
\newcommand{\point}{\,\cdot\,}

\renewcommand{\a}{\alpha}

\newcommand{\e}{\text{e}}

\newcommand{\Tree}{\mathbb T}

\newcommand{\Tor}{T_{\text{or}}}
\newcommand{\Pnt}{\mathbb P_n^{t}}
\newcommand{\Hetoile}{H^\diamond} 
\newcommand{\Pni}{\P_n^{(i)}}
\newcommand{\Pninf}{\P_n^{(\infty)}}
\newcommand{\hni}{h_n^{(i)}}
\newcommand{\hi}{h^{(i)}}
\newcommand{\cvd}{\Rightarrow}
\newcommand{\Loi}[2]{\mathcal L(#1,\; #2)}
\newcommand{\egalLoi}{\stackrel{\mathcal L}{=}}
\newcommand{\pii}{\pi^{(i)}}
\newcommand{\Tnk}[2]{T_{#1,#2}}
\newcommand{\Ti}[1]{T^{(i)}_{#1}}
\newcommand{\hTnk}[2]{\hat T_{#1,#2}}
\newcommand{\Eni}{\E_n^{(i)}}
\newcommand{\Eninf}{\E_n^{(\infty)}}
\newcommand{\phii}[2]{\phi_{#1,#2}}
\newcommand{\psii}[2]{\psi_{#1,#2}}
\newcommand{\Tori}{\Tor^{(i)}}
\renewcommand{\sp}{p}

    \makeatletter
    \let\@fnsymbol\@arabic
    \makeatother

\begin{document}
\selectlanguage{english}
\setlength{\parindent}{0cm}
\title{Sample genealogy and mutational patterns for critical branching populations}
\author{Guillaume Achaz\footnotemark[1] ${}^,$\footnotemark[2] ${}^,$\footnotemark[3] \and Cécile Delaporte\footnotemark[3] ${}^,$\footnotemark[4] \and Amaury Lambert\footnotemark[3] ${}^,$\footnotemark[4]}

\date{}
\maketitle

\footnotetext[1]{UMR 7138, Evolution Paris-Seine, UPMC \& CNRS, Paris.}
\footnotetext[2]{Atelier de Bioinformatique, UPMC, Paris.}
\footnotetext[3]{UMR 7241, Centre Interdisciplinaire de Recherche en Biologie, Collège de France, Paris.}
\footnotetext[4]{UMR 7599, Laboratoire de Probabilités et Modèles Aléatoires, UPMC \& CNRS, Paris.}
\renewcommand{\thefootnote}{\*} 
\footnote{{Corresponding author :} amaury.lambert@upmc.fr}

\begin{abstract}
 We study a universal object for the genealogy of a sample in populations with mutations: the critical birth-death process with Poissonian mutations, conditioned on its population size at a fixed time horizon. We show how this process arises as the law of the genealogy of a sample in a large class of critical branching populations with mutations at birth, namely populations converging, in a large population asymptotic, towards the continuum random tree. We extend this model to populations with random foundation times, with (potentially improper) prior distributions $g_i:x\mapsto x^{-i}$, $i\in\Z_+$, including the so-called \textit{uniform} ($i=0$) and \textit{log-uniform} ($i=1$) priors.\\

 We first investigate the mutational patterns arising from these models, by studying the site frequency spectrum of a sample with fixed size, i.e. the number of mutations carried by $k$ individuals in the sample. Explicit formulae for the expected frequency spectrum of a sample are provided, in the cases of a fixed foundation time, and of a uniform and log-uniform prior on the foundation time. Second, we establish the convergence in distribution, for large sample sizes, of the (suitably renormalized) tree spanned by the sample genealogy with prior $g_i$ on the time of origin. We finally prove that the limiting genealogies with different priors can all be embedded in the same realization of a given Poisson point measure.
\end{abstract}
\par\bigskip

\noindent\textit{2010 AMS Classification : }92D10, 60J80 (Primary), 92D25, 60F17, 60G55, 60G57, 60J85 (Secondary)\\

\section{Introduction}\label{sec_new_intro}

\qquad A major concern in population genetics is the prediction of patterns of genetic variation with help of stochastic models. The reference model currently used by biologists to answer this question is the Kingman coalescent model \cite{Kingman_a,Kingman_b} coupled with Poissonian mutations on the lineages. As the scaling limit of numerous constant population size models, such as Wright-Fisher and Moran models, it encompasses the two population models that are most commonly used by biologists. The genealogical structure of a sample (rather than of the total population) is well-known (equivalently given by the Kingman coalescent), and explicit results on the allelic partition generated by rare, neutral mutations (equivalent to a Kingman coalescent with Poissonian mutations) are provided by Ewens' sampling formula \cite{Ewens,Durrett}. In this work, we intend to study the genealogical and mutational patterns of a sample from a branching population, in order to offer an alternative model where the constant 
population size assumption is released, with no \textit{a priori} assumption on the variation of the population size 
over time. The sampling is here essential to make the model applicable to real data and comparable to the Kingman coalescent model.\\ 

\qquad The genealogy of branching populations was in particular studied by L. Popovic in \cite{Popovic}, in the setting of the critical birth-death process conditioned on its population size at a fixed time horizon, and later by A. Lambert in \cite{ALContour} in the more general framework of splitting trees. The genealogy of the extant individuals is described as a random point process, called \textit{coalescent point process}, which distribution is characterized by a sequence of i.i.d. random variables. \\

\qquad Here we want to focus on the genealogy of a sample rather than of the total extant population. 
The question of sampling in birth-death models has already been approached with two different points of view. On the one hand, \cite{Popovic} and \cite{StadlerIncomplete} deal with Bernoulli sampling of the total population. This approach rather applies to the species scale, for example in the case of incomplete phylogenies. On the other hand, in \cite{StadlerLTT} and \cite{StadlerIncomplete}, T. Stadler considers the case of a uniform sample of $m$ individuals among the extant ones, in the birth-death process conditioned on its population size at present time, with uniform prior on its time of origin. Our approach is based on Bernoulli sampling with conditioning on the sample size, in order to get a uniform sample with fixed size without having to condition on the total extant population size.\\

We first consider in Section \ref{sec_new_intro_1} sample genealogies in a general framework of branching populations with neutral mutations at birth. We make use of convergence results obtained by one of the authors \cite{LPWM2} to show how a broad class of such populations all result in the same distribution for the genealogy of a sample, namely the law of a critical birth-death model with Poissonian mutations on the lineages. We then specify in Section \ref{sec_new_intro_2} the model that we adopt for the rest of the paper. We finally present in Section \ref{sec_new_intro_3} the outline and the main results of this work : in Section \ref{sec_1}, we investigate the law of the genealogy of a sample in the critical birth-death model conditioned on its sample size, with various prior distributions on the foundation time of the population. We provide in Section \ref{sec_2} explicit formulae for the expected site frequency spectrum of the sample. Section \ref{sec_3} is then devoted to the convergence in 
distribution of the sample genealogy, as the sample size gets large. Furthermore, we state that the limiting genealogies with different priors can all be embedded in the same realization of a given Poisson point measure.

\subsection{Genealogies and sampling in branching populations conditioned on survival}\label{sec_new_intro_1}


\qquad Let us first consider a very general model of branching populations with mutations : let ($\Tree_N)_{N\in\N}$ be a sequence of splitting trees, i.e. random trees where individuals have lifetimes that do not necessarily follow an exponential distribution, during which they give birth at constant rate to i.i.d copies of themselves \cite{Geiger,GK,ALContour}. For any $N$, $\Tree_N$ is characterized by its so-called \textit{lifespan measure} $\L_N$, which is a $\sigma$-finite measure on $(0,\infty)$ such that $\int (1\wedge r)\L_N(\d r)<\infty$. We further assume that any individual in $\Tree_N$ experiences, conditional on her lifetime $r$, a mutation at birth with probability $f_N(r)$, where $f_N$ is a continuous function from $\R_+^*$ to $[0,1]$ called \textit{mutation function}. We adopt the classical assumptions of neutral mutations (i.e. mutations do not affect the population dynamics) and of the infinite-site model \cite{KimuraISM} : each individual is associated to a DNA sequence, and each mutation 
occurs 
at a site that has never mutated before. \\

\begin{figure}[!h]
 \begin{center}
  \includegraphics[width=.98\linewidth]{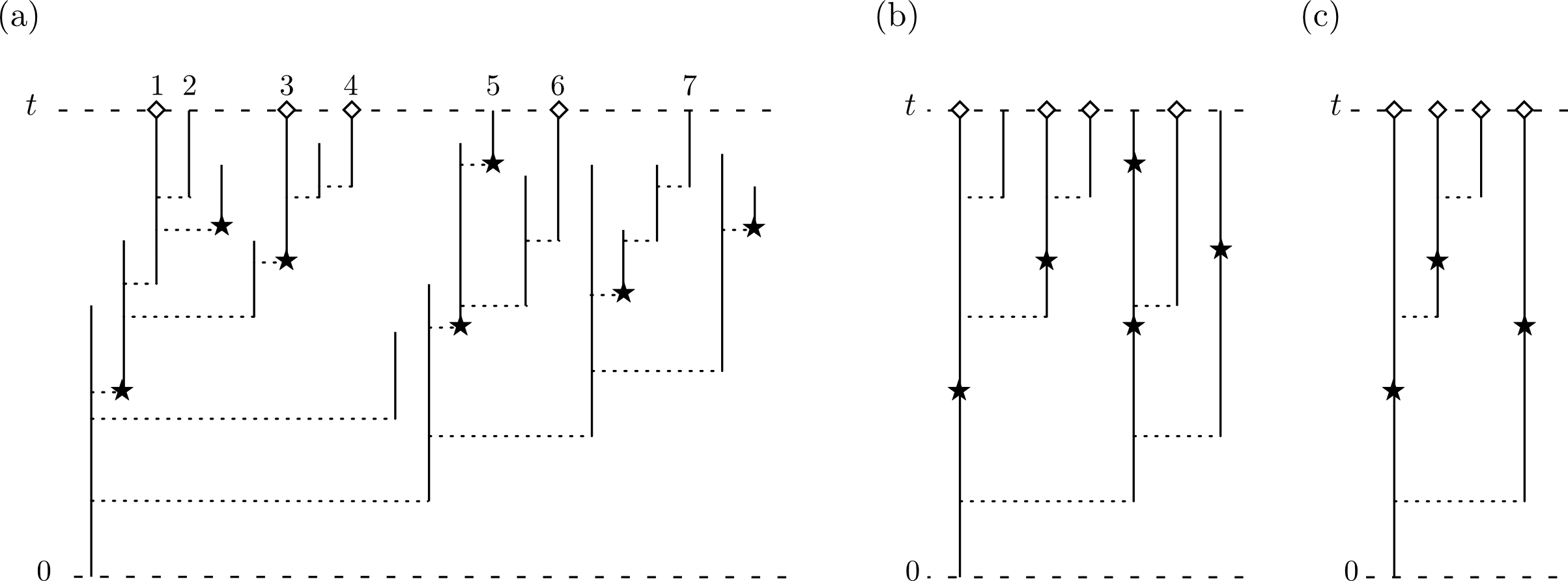}
  \caption{In the three panels (a), (b), (c), the vertical axis indicates time. The horizontal (dotted) lines show filiation. Mutations are symbolized by $\star$ and sampled individuals by $\diamond$. \\
(a) An example of the rescaled tree $\Tree_N$ with $7$ extant individuals at $t$, where $4$ individuals are sampled. \\
(b) its (marked) coalescent point process (later referred to as $\Sigma_N$),\\
(c) and the (marked) coalescent point process of the sampled individuals.}
  \label{fig_new_arbre_cpp_ech}
 \end{center}
\end{figure}

Finally, we fix $t>0$, and we condition $\Tree_N$ on survival at time $Nt$. We work later in a time scale where a unit of time is proportional to $N$ : the factor $N$ can thus be seen as a counterpart of the constant population size of the Wright-Fisher model. We assume that each individual alive at $Nt$ is independently sampled with probability $\sp_N\in(0,1)$. Individuals are labeled according to the order defined in \cite[Sec. 1.1]{ALallelicpartition} (\og left to right\fg\ order associated to the planar representation of the tree when daughters all sprout to the right of their mother), 
and we denote by $I_N=(I_{Nj})_j$ the sequence of indices of the sampled individuals. See Figure \ref{fig_new_arbre_cpp_ech} for a graphical representation of $\Tree_N$, and of some objects hereafter defined.\

\qquad We are here interested in the distribution of the genealogy of the sampled individuals in $\Tree_N$, and we consider the model under two slightly different points of view : in case (I), relying on results of \cite{LPWM2}, we consider a scaling limit in a large population asymptotic, while in case (II), we consider the example of the critical birth-death process, for which results can be obtained without necessarily having to consider $N\to\infty$. We show here how these two settings lead to the same distribution for the genealogy of a sample, justifying hence the model we later consider for the rest of the paper.\\

To this aim we rescale time in $\Tree_N$ by multipying all the edge lengths of $\Tree_N$ by a factor $1/N$. This rescaled tree is still denoted by $\Tree_N$, and is now originating at time $t$. Then we introduce, for any $N\in\N$, the so called \textit{marked coalescent point process} $\Sigma_N$ \cite{LPWM2},  i.e. the tree spanned by the genealogy of the extant population of $\Tree_N$ at time $t$, enriched with the mutational history of extant individuals. More precisely, $\Sigma_N$ is a point measure that can be expressed as $\Sigma_N=\sum_{i=1}^{\mathcal N-1}\delta_{(i,\sigma_i^N)}$ where $\mathcal N$ is the number of extant individuals at time $t$, and for any $1\leq i\leq \mathcal N-1$, $\sigma_i^N$ is itself a point measure, whose set of atoms contains, in addition to the coalescence time between individuals $i$ and $i+1$, all the times at which a mutation occurred on the $i$-th lineage (see Figure \ref{fig_new_arbre_cpp_ech}).\\ 

\paragraph{(I) Scaling limit.} \quad \\
First, we assume that $(\Tree_N)$ converges, as $N\to\infty$, towards a Brownian tree (see e.g. \cite{Aldous}) : for any $N\in\N$, for any $\lambda\geq0$, define $\psi_N(\lambda):=-\lambda-\int_{(0,\infty)}(1-e^{-\lambda r})\L_N(\d r)$. 
We assume that the sequence $(\Tree_N)$ follows (a particular case of) Assumption \HypLevy\ in \cite{LPWM2} : 
\begin{center}
\begin{minipage}{.9\linewidth}
\textit{Assumption \HypLevy\ : There exists a sequence of positive real numbers $(d_N)_{N\geq1}$ such that as $N\to\infty$, the sequence  $(d_N\psi_N(\point/N))$ converges towards $\lambda\mapsto\lambda^2$, i.e. the Laplace exponent of a Brownian motion.}
\end{minipage} 
\end{center}
This assumption has to be interpreted as the convergence in law of the so-called \textit{jumping chronological contour process} of the rescaled tree $\Tree_N$, which distribution is characterized by a Lévy process with finite variation, drift $-1$ and Lévy measure $\L_N$ \cite{ALContour}. 

 Second, we fix $\theta\in\R_+$ and we suppose that the sequence of mutation functions $(f_N)$ satisfies one of the following convergence assumptions \cite{LPWM2} : 

\begin{center}
\begin{minipage}{.9\linewidth}
\textit{Assumption \HypMutConst\ : For all $N\geq1$, for all $u\in\R_+$, $f_N(u)=\theta_N$, where $\theta_N\in[0,1]$ is such that $\frac{d_N}{N}\theta_N \underset{N\to\infty}\longrightarrow\theta$.}
\end{minipage} 
\end{center}

\begin{center}
\begin{minipage}{.9\linewidth}
\textit{Assumption \HypMut\ : The sequence $\big(u\mapsto \frac{f_N(Nu)}{1\wedge u}\big)$ converges uniformly to $u\mapsto \frac{f(u)}{1\wedge u}$ on $\R_+^*$, where $f$ is a continuous function from $\R_+$ to $\R_+$ satisfying $f(u)/u \to \theta$ as $u\to0^+$.}
\end{minipage} 
\end{center}

 Then we have the following convergence.
\begin{theoreme*}{\cite[Th.3.2]{LPWM2}} 
 The (space rescaled) point measure $\Sigma_N=\sum_{i=1}^{\mathcal N-1}\delta_{(i\frac{d_N}N,\sigma_i^N)}$ converges in distribution, as $N\to\infty$, towards a Poisson point process on \upshape$[0,\e]\times (0,t)$ \itshape with intensity \upshape$\d l\,x^{-2}\d x$\itshape, where \upshape$\e$ \itshape is an independent exponential variable with parameter $1/t$, with independent Poissonian mutations at rate $\theta$ on the lineages.
\end{theoreme*}
Besides, we assume that the sampling probability is given by $\sp_N=\sp\, N/d_N$, where $\sp$ is a fixed positive real number such that $\sp_N\in(0,1)$ for $N$ large enough. Then the rescaled sequence $(\frac {d_N}N I_N)$ of indices of the sampled individuals (independent of $(\Tree_N)$), converges towards the sequence of jump times of an independent Poisson process with rate $\sp$. The joint convergence of $\Sigma_N$ with $\frac {d_N}N I_N$ is of course provided by their independence. 
\par\medskip
As a consequence, from \cite{ALallelicpartition} we deduce that the coalescent point process of the sampled individuals is then distributed as the coalescent point process of a critical birth-death model with rate $\sp$ conditioned on survival at time $t$, with independent Poissonian mutations at rate $\theta$ on the lineages. 

\paragraph{(II) Critical birth-death tree.} \quad \\
Second, fix $N\in\N$, $p\in(0,N)$, and consider the example where $\Tree_N$ is a critical birth-death tree with rate $N$ conditioned on survival at time $t$. Then, set $p_N=p/N$ and assume that the mutation function $f_N$ is constant, equal to $\theta/N$. This is in fact a particular case of (I) (Assumptions \HypLevy\ and \HypMutConst\ are satisfied with $d_N=N^2$), but here we do not need to let $N\to\infty$. For any $N\in\N$, the marked coalescent point process $\Sigma_N$ is distributed as the coalescent point process of a critical birth-death model with rate $1$ conditioned on survival at time $t$, with Poissonian mutations at rate $\theta$ on the lineages (see \cite[Sec.3]{Popovic} and \cite[Ex.1]{LPWM2}). Finally, from \cite{ALallelicpartition}, we get that the coalescent point process of the sample is then distributed, exactly as above, as the coalescent point process of a critical birth-death model with rate $\sp$ conditioned on survival at time $t$, with independent Poissonian mutations at rate $\
theta$ on the lineages.\\

\par\bigskip

\qquad Since the two cases (I) and (II) result in the same distribution for the genealogy of a sample, we limit our study to case (II). Besides, since the mutation schemes arise as independent of genealogies, the results concerning distributions of genealogies are stated without reference to mutations. \\

\subsection{Model with conditioning on the sample size}\label{sec_new_intro_2}

\qquad From now on, consider $\Tree$ a critical birth-death tree with rate $1$. Time is now counted backwards into the past, i.e. \og present time\fg\ is now time $0$, and \og$u$ units of time before present\fg\ is now time $u$. We begin with the case of a fixed foundation time of the population. The model has four parameters : a time $t\in\R_+^*$, a scaling factor $N\in\R_+^*$, a positive integer $n$ (the sample size), and a sampling parameter $\sp\in(0,N)$.\\

Assume first that $\Tree$ has been founded $Nt$ units of time ago. As previously, individuals are independently sampled at present time, with probability $\sp/N$. Besides, we rescale time by a factor $1/N$ (all the edge lengths are then multiplied by a factor $1/N$). We keep the notation $\Tree$ for the rescaled tree, so that $\Tree$ is now a critical birth-death tree with rate $N$, originating at time $t$. \\

\begin{figure}[!h]
 \begin{center}
  \includegraphics[width=.8\linewidth]{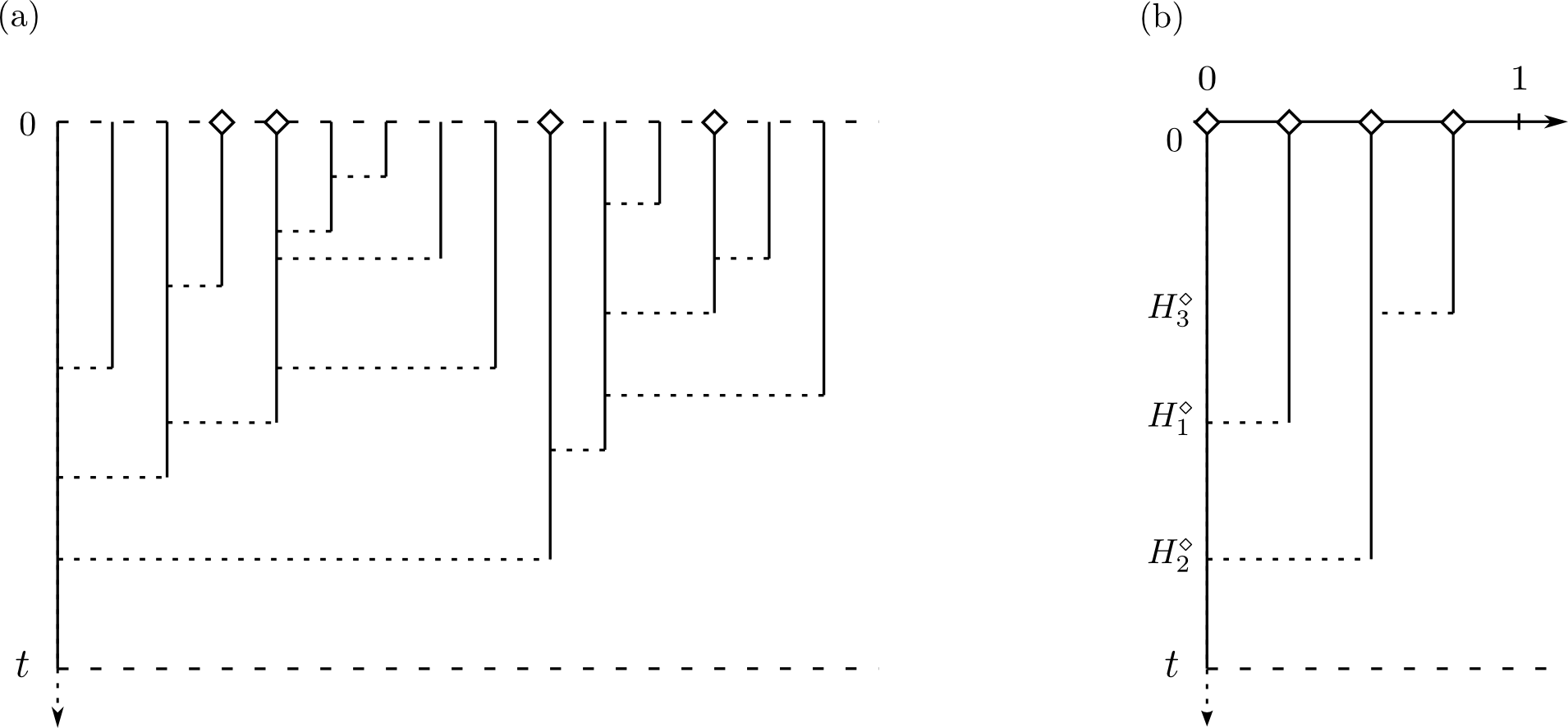}
  \caption{In both figures (a) and (b), the vertical axis indicates time (running backwards). \\
(a) A graphical representation of the coalescent point process at present time of a (rescaled) tree $\Tree$ originating at time $t$ with $15$ extant individuals and $n=4$ sampled individuals (symbolized by $\diamond$). The horizontal lines show filiation.\\
(b) A graphical representation of the coalescent point process $\pi_n=\sum_{k=1}^{n-1} \delta_{\left(\frac kn,\Hetoile_k\right)}$ of the sample represented in (a).}
  \label{fig_new_cpp_ech}
 \end{center}
\end{figure}

We now introduce the conditioning on the sample size : we condition $\Tree$ on having $n$ sampled individuals at present time. Note that after conditioning, the distribution of the $n$ sampled individuals within the total extant population does not depend on $\sp$, and is a posteriori equivalent to uniform, sequential sampling. \\

The genealogy of the $n$ sampled individuals is characterized by its coalescent point process 
$$\pi_n=\sum_{k=1}^{n-1} \delta_{\left(\frac kn,\Hetoile_k\right)},$$
where for $1\leq k\leq n-1$, $\Hetoile_k$ is the divergence time between the $k$-th and the $(k+1)$-th sampled individual in the rescaled tree $\Tree$ (see Figure \ref{fig_new_cpp_ech}). The space rescaling by a factor $1/n$ ensures in particular that the supports of the measures $\pi_n$ converge as $n\to\infty$, which is required by the results later established in the large sample size asymptotic. Besides, recall that thanks to their independence with the genealogy, mutations are for now deliberately omitted. Finally, we define $(\Tnk nk)_{1\leq k\leq n-1}$ the decreasing reordering of the divergence times $(\Hetoile_k)_{1\leq k\leq n-1}$.

\subsection{Outline and statement of results}\label{sec_new_intro_3}

\qquad The first purpose of this paper is to study the distribution of $\pi_n$, under various hypotheses on the origin of the process : we denote by 
\begin{itemize}
 \item $\Pnt$ the law of the rescaled tree $\Tree$ with fixed time of origin $t$ and sample size $n$,
 \item $\Pninf$ the law of $\Tree$ with infinite time of origin and sample size $n$,
 \item $\Pni$ the law of $\Tree$ with random time of origin, with (potentially improper) prior distribution $g_i:x\mapsto x^{-i}$, $i\in\Z_+$, and sample size $n$.
\end{itemize}
Note that the case $i=0$ corresponds to the case of a uniform prior investigated in \cite{Aldous-Popovic} and \cite{StadlerNew}. This study, presented in Section \ref{sec_genealogy}, will then enable us to derive results concerning mutational patterns of the sample (Section \ref{sec_mut}), and then concerning the behaviour of the genealogy as the sample size gets large (Section \ref{sec_cv}). \\

\subsubsection{A universal law for the genealogy of a sample}\label{sec_1}

\qquad First, in the case of a fixed time of origin, the law of $\pi_n$ under $\Pnt$ is independent of $N$ and is specified by the following result (Theorem \ref{th_loi_CPP_Tor_fixe}) :
\begin{theoreme*}
 Under $\Pnt$, $(\Hetoile_i)_{1\leq i\leq n-1}$ is a sequence of i.i.d. random variables with probability density function $x\mapsto \frac\sp{(1+\sp x)^2}\frac{1+\sp t}{\sp t}\mathbbm1_{(0,t)}(x)$. In other words, the coalescent point process $\pi_n$ has the law of the genealogy of a critical birth-death tree with rate $\sp$ conditioned on having $n$ extant individuals at time $t$.
\end{theoreme*}
We then prove that this equality in law still holds when letting the time go to infinity or when randomizing the time of origin (with prior distribution $g_i$, $i\in\Z_+$) in both processes : for example, under $\Pni$ the coalescent point process $\pi_n$ has the law of the genealogy of a critical birth-death tree with rate $\sp$, with prior $g_i$ on its time of origin, and conditioned on having $n$ extant individuals at present time. Hence whatever the assumption on the foundation time of the population, the study of the genealogy of the sample boils down to the same object : the genealogy of a critical birth-death process with rate $\sp$, with extant population size $n$.

Following on from results provided by \cite{StadlerNew} in the case of a uniform prior, we then obtain the following property for the successive divergence times $(\Tnk nk)_{1\leq k\leq n-1}$ (Proposition \ref{prop_moments_Eni(Tnk)}) :
\begin{proposition*}
 Under $\Pni$, the time $\Tnk nk$ to the $k$-th most recent common ancestor 
has finite moment of order $m$ iff $m\leq k+i$.
\end{proposition*}
\qquad Although we limited here our study to the framework (II) introduced earlier, one could certainly generalize these results (and the upcoming ones) to the scaling limit of case (I). To prove this, one would have to consider a sequence of trees conditioned on their sample size, and then to establish the convergence, in the large population asymptotic, of the marked coalescent point process of the sample. This is however beyond the scope of the present paper. \\

\subsubsection{Mutational patterns}\label{sec_2}
\qquad In Section \ref{sec_mut}, we study the so-called \textit{site frequency spectrum} of the sample, i.e. the $(n-1)$-tuple $(\xi_1,\ldots,\xi_{n-1})$, where $\xi_k$ is the number of mutations carried by $k$ individuals in the sample.  Various results for the frequency spectrum in the framework of general branching processes are established in \cite{ALallelicpartition,CL1,CL2,Richard}. One of the authors investigates in \cite{ALallelicpartition} the case of coalescent point processes with Poissonian mutations on germ lines, and obtains asymptotic results for the site and allele frequency spectrum of large samples. Explicit formulae for the expected allele frequency spectrum of a splitting tree with $n$ individuals at fixed time horizon $t$ are provided by N. Champagnat and this author in the case of Poissonian mutations on the lineages \cite{CL1}, and by M. Richard in the case of mutations at birth \cite{Richard}. Their results are compared in \cite{CLR} in the particular case of birth-death processes. 
Further 
results about the asymptotic behaviour, as $t\to\infty$, of large (resp. old) families, i.e. families with most frequent (resp. oldest) types, are developed in \cite{CL2}.\\

In this article, we get explicit formulae for the expected site frequency spectrum $(\xi_k)_{1\leq k\leq n-1}$ of the sample under $\Pnt$, $\Pninf$, $\P_n^{(0)}$ and $\P_n^{(1)}$. According to Section \ref{sec_new_intro_1}, mutations are assumed to occur at constant rate $\theta$ on the lineages. Two different methods are used to obtain the expectation of the $\xi_k$. On the one hand, the similarity of the model with \cite{ALallelicpartition} allows us to make use of a proof method developed in this article. Indeed, according to the results of Section \ref{sec_genealogy}, the framework used in \cite{ALallelicpartition} covers our setting in the case of an infinite time of origin. On the other hand, for each $k$, $\E(\xi_k)$ can be expressed as a linear combination of the expectations of branching times \cite{Wakeley}. Although the first method could be used to prove all the results of this section, the second one provides very short proofs in the cases of an infinite time of origin and of a uniform prior. 
First 
under $\Pninf$, the absence of a first moment for the time to the most recent common ancestor yields immediately the 
following result (Proposition \ref{prop_Eninfty(xi_k)}).
\begin{proposition*}
For any $k\in\{1,\ldots,n-1\}$, $\E_n^{(\infty)}(\xi_k)$ is infinite.
\end{proposition*}
Second, using the fact that the expected divergence times, under the Kingman coalescent model, and under the (suitably rescaled) critical birth-death process with uniform prior on its time of origin, are equal \cite{StadlerNew}, we deduce that the expected site frequency spectrum under $\P_n^{(0)}$ is that of a sample of the Kingman coalescent \cite[(4.20)]{Wakeley} (Proposition \ref{prop_En0(xi_k)}).
\begin{proposition*}
 For any $k\in\{1,\ldots,n-1\}$, $\E_n^{(0)}(\xi_k)=n\theta/k\sp$. 
\end{proposition*}
Finally, the formulas obtained in the remaining two cases are the following (Propositions \ref{prop_Ent(xi_k)} and \ref{prop_En1(xi_k)}).
\begin{proposition*}
 For any $k\in\{1,\ldots,n-1\}$, $t\in\R_+^*$, defining $\tau:=\sp t$, we have
\begin{multline*}
 \E_n^t(\xi_k)
=\frac\theta\sp \ \Bigg\{\frac{n-3k-1}{k}+\frac{(n-k-1)(k+1)}{k\tau} \\
 \quad+\frac{(1+\tau)^{k-1}}{\tau^{k+1}}\Big[2\tau^2-(n-2k-1)2\tau-(n-k-1)(k+1)\Big]\bigg[\ln(1+\tau)-\sum_{i=1}^{k-1}\frac1i\left(\frac{\tau}{1+\tau}\right)^{i}\bigg] \Bigg\}.
\end{multline*}
\end{proposition*}

\begin{proposition*}
  For any $k\in\{1,\ldots,n-3\}$, 
$$\E_n^{(1)}(\xi_k)=\frac\theta\sp\,\frac{n(n-1)}{(n-k)(n-k-2)}\Bigg[\frac{n+k-2}{k}-\frac{2(n-1)}{n-k-1}(\mathcal H_{n-1}-\mathcal H_k)\Bigg],$$
where for any $k\in\N$, $\mathcal H_k=\sum_{j=1}^k j^{-1}$.
\end{proposition*}

\subsubsection{Convergence of genealogies for large sample sizes}\label{sec_3}

\qquad We investigate in Section \ref{sec_cv} the asymptotic behaviour of the coalescent point process $\pi_n$, as $n\to\infty$. We take inspiration from asymptotic results presented in \cite{Popovic} and \cite{Aldous-Popovic}. First, L. Popovic obtains in \cite{Popovic} the convergence of the (suitably rescaled) coalescent point process of a critical birth-death process conditioned on its population size at time $t$ towards a certain Poisson point measure on $(0,1)\times(0,t)$. 
Using this result, she then obtains with D. Aldous in \cite{Aldous-Popovic} a similar convergence for the model with uniform prior on the time of origin. Here we extend this to the cases of an infinite time of origin, and of a random time of origin with prior $g_i$, $i\in\N$.\\

Obtaining such asymptotic results requires to let the sampling parameter $\sp$ depend on $n$ in such a way that $\sp=n/\a$, with $\a>0$. It ensures indeed that the expected number of sampled individuals is of the order of the sample size $n$. We then obtain the following convergences (Theorem \ref{th_cv_pin}).
\begin{theoreme*} Denote by $\pi^t$ the Poisson point measure with intensity \upshape$\a\d l\, x^{-2}\d x \mathbbm1_{(l,x)\in(0,1)\times(0,\a t)}$\itshape. 
\begin{enumerate}[\upshape a)\itshape ]
 \item Under $\Pninf$, the coalescent point process $\pi_n$ converges in law, as $n\to\infty$, towards the Poisson point measure $\pi$ with intensity measure \upshape$\a\d l\,x^{-2}\d x$ \itshape on $(0,1)\times\R_+^*$.
 \item For any $i\in\Z_+$, under $\Pni$, the joint law of the time of origin, along with $\pi_n$, converges  as $n\to\infty$ towards a pair \upshape$(\Tori,\,\pii)$\itshape, such that \upshape$\Tori$ \itshape follows an inverse-gamma distribution with parameters $(i+1,\a)$, and conditional on \upshape$\Tori=t$\itshape, $\pii$ is distributed as $\pi^t$.
\end{enumerate} 
\end{theoreme*}

The last result we obtain describes the links between the different random measures obtained in the limit. Let us order the atoms of our point processes w.r.t. their second coordinate. We prove that the random variable $\Tori$ is distributed as the $(i+1)$-th largest atom of the Poisson point process $\pi$, and we then deduce the following theorem (Theorem \ref{th_pi_plus_gd_atome}).
\begin{theoreme*}
 The point measure $\pii$ is distributed as the point process obtained from $\pi$ by removing its $i+1$ largest atoms. 
\end{theoreme*}
 In other words, genealogies with different priors can all be embedded in the same realization of the point measure $\pi$.
\par\bigskip
\section{A universal distribution for the genealogy of a sample}\label{sec_genealogy}

Let us consider the model defined in Section \ref{sec_new_intro_2} and specify some notation. Recall that the rescaled tree $\Tree$ is a critical birth-death tree with parameter $N$ originating at time $t$, and that each extant individual in $\Tree$ is independently sampled with probability $\sp/N$. \\

We denote by $\mathcal N$ the number of extant individuals at present time in $\Tree$, and we label these individuals from $1$ to $\mathcal N$, using the order defined in \cite[Sec. 1.1]{ALallelicpartition}. In order to formalize the sampling process, we introduce a sequence $(I_j)_{j\geq1} $ of random variables, such that $(I_1,I_2-I_1,I_3-I_2,...)$ forms a sequence of i.i.d. geometric random variables with success probability $\sp/N$. Then for any $j$ such that $I_j\leq\mathcal N$, $I_j$ is the label of the $j$-th sampled individual in the extant population at present time (in the previously defined order). The conditioning on the sample size to be equal to $n$ means thus conditioning on $\{I_n\leq \mathcal N<I_{n+1}\}$.\\

\begin{figure}[!h]
\begin{center}
  \includegraphics[width=0.98\linewidth]{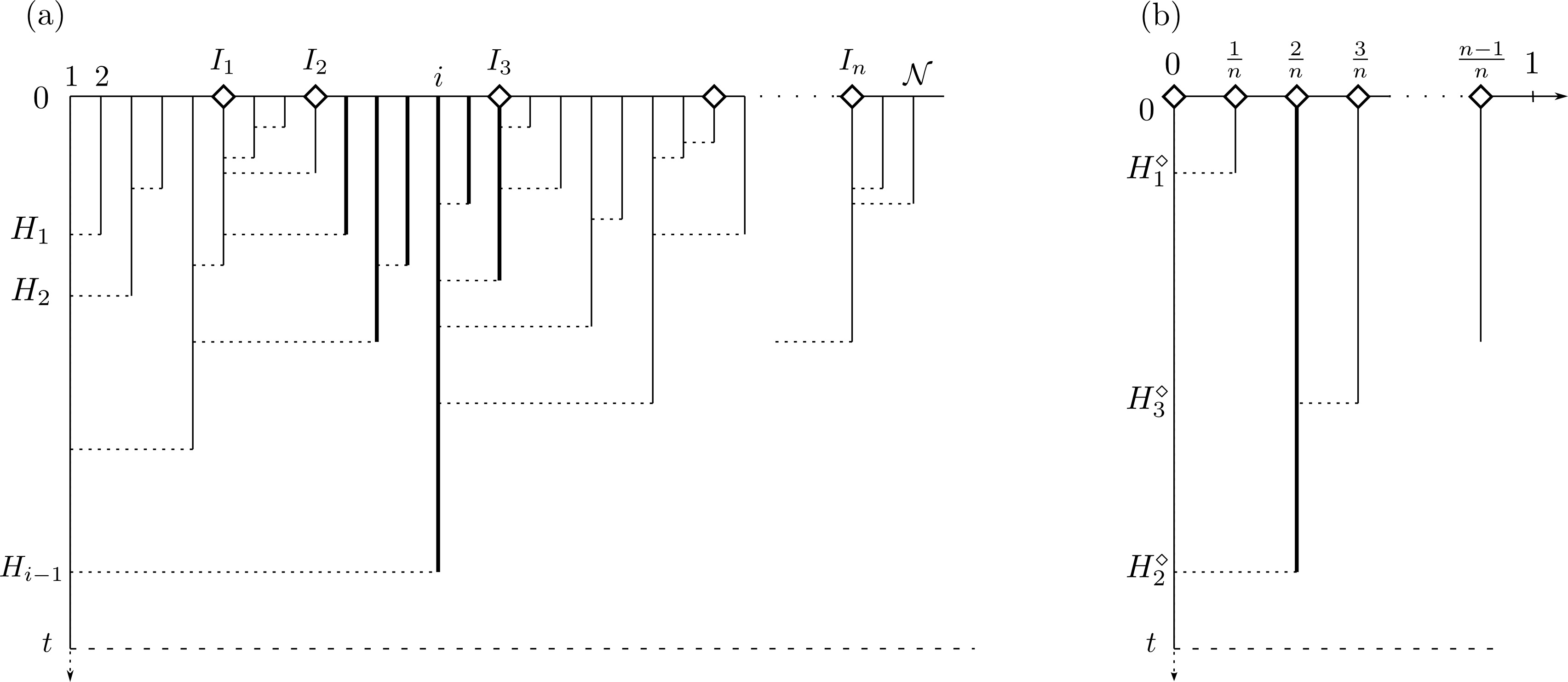}
  \caption{
(a)The coalescent point process at present time of a (rescaled) population originating at time $t$ with $n$ sampled individuals (symbolized by $\diamond$). The $\mathcal N$ vertical branches represent the sequence $(H_i)_{1\leq i\leq\mathcal N}$.\\
(b) The coalescent point process $\pi_n$ of the sample represented in figure (a). The equality $\Hetoile_2=\max\{H_{I_2+1},\ldots,H_{I_3}\}$ is illustrated by bold lines.}
   \label{fig_model}
\end{center}
\end{figure}

Let us now explain the link between the genealogy of the total extant population and the genealogy of the sample. Denote by $(H_i)_{1\leq i\leq \mathcal N-1}$ the sequence of node depths of the coalescent point process of the total extant population, i.e. for any $1\leq i\leq \mathcal N-1$, $H_i$ is the divergence time between individual $i$ and individual $i+1$ in the rescaled tree $\Tree$. We know from \cite[Th.5.4]{ALContour} that for any $1\leq i<j\leq \mathcal N$, the divergence time between individual $i$ and $j$ is given by the maximum of the node depths $\{H_{i+1},\ldots,H_j\}$. As a consequence, the divergence time $\Hetoile_i$ between individual $I_i$ and individual $I_{i+1}$ in $\Tree$, $1\leq i\leq n-1$, is given by 
$$\Hetoile_i=\max\{{H_{I_i+1}},\ldots,{H_{I_{i+1}}}\}.$$

Finally we recall the definition of the point measure $\pi_n$ : 
$$\pi_n=\sum_{k=1}^{n-1} \delta_{\left(\frac kn,\Hetoile_k\right)}.$$
In the sequel we equally call \og coalescent point process\fg\ of the sample, the measure $\pi_n$ and the sequence $(\Hetoile_i)_{1\leq i\leq n-1}$. See Figure \ref{fig_model} for a graphical representation of the objects defined above.\\

The aim of this section is to characterize the law of the genealogy of the sample, under different assumptions on the time of origin. Section \ref{sec_Tor_fixe} establishes the distribution of $\pi_n$ in the case of a fixed (possibly infinite) time of origin. In Section \ref{sec_Tor_prior}, we randomize the time of origin by giving it a prior distribution of the form $x\mapsto x^{-i}$, $i\in\Z_+$.

\subsection{Fixed time of origin}\label{sec_Tor_fixe}

We denote by $\P^t$ the law of the rescaled tree $\Tree$ originating at time $t$, and we recall that $\Pnt$ denotes the law of $\Tree$ originating at time $t$ and conditioned on having $n$ sampled individuals at present time, i.e on $\{I_n\leq \mathcal N<I_{n+1}\}$. The following theorem specifies the law of the sample genealogy under $\Pnt$.

\begin{theoreme}\label{th_loi_CPP_Tor_fixe}
Under $\Pnt$, the coalescent point process $(\Hetoile_i)_{1\leq i\leq n-1}$ is a sequence of i.i.d. random variables with probability density function $$x\mapsto \frac\sp{(1+\sp x)^2}\frac{1+\sp t}{\sp t}\mathbbm1_{(0,t)}(x).$$ 
\end{theoreme}

\begin{remarque}
According to \cite[Lem.3]{Popovic}, the rescaled coalescent point process of the $n$ sampled individuals is thus distributed as the coalescent point process of the population at time $t$ of a critical branching process with rate $\sp$, conditioned on having $n$ extant individuals at time $t$ -- or equivalently, as the coalescent point process of the population at time $\sp t$ of a critical branching process with rate $1$, conditioned on having $n$ extant individuals at time $\sp t$, and then rescaled by a factor $1/\sp$. 
\end{remarque}
\begin{remarque}\label{rem_effet_ptN}
 It is interesting to note that the independence w.r.t. $N$ of the law of $\pi_n$ under $\Pnt$ implies that the parameter $N$ has only a scaling effect on the law of the genealogy. On the contrary, the parameters $\sp$ and $t$ both affect the branch lengths ratios, through the conditioning on the population size at a fixed time.
\end{remarque}

We extend the theorem to the limiting case $t\to\infty$ : recall that $\Pninf(\Tree\in\point)=\underset{t\to\infty}\lim \Pnt(\Tree\in\point)$. We have the following statement.
\begin{proposition}\label{prop_loi_CPP_Tor_infini}
 Under $\Pninf$, $(\Hetoile_i)_{1\leq i\leq n-1}$ is a sequence of i.i.d. random variables with density function $x\mapsto \frac\sp{(1+\sp x)^2}\mathbbm1_{\R_+}(x)$.
\end{proposition}
\par\bigskip

Recall that for any $1\leq k\leq n-1$, $\Tnk nk$ is defined as the $k$-th order statistic of the sequence $(\Hetoile_i)_{1\leq i\leq n-1}$. In particular, $\Tnk n1$ is the time to the most recent common ancestor of the sample. The following proposition provides the $m$-th moment of $\Tnk nk$ under $\Pninf$.

\begin{proposition}\label{prop_moments_Eninfty(Tnk)}
 For any $1\leq k\leq n-1$ and $m\geq1$, the $m$-th moment of $\Tnk nk$ under $\Pninf$ is finite iff $m\leq k-1$. Specifically, for $m\leq k-1$,
 $$\E_n^{(\infty)}((\Tnk nk) ^m)=\frac{\binom{n-k+m-1}m}{\sp^m \binom{k-1}m}.$$
\end{proposition}

In particular, the time to the most recent common ancestor has infinite expectation under $\Pninf$.\\

\begin{demopr}{\ref{prop_moments_Eninfty(Tnk)}}
Using the definition of  $\Tnk nk$ as the $k$-th order statistic of the i.i.d. random variables $(\Hetoile_i)_{1\leq i\leq n-1}$ with density function $x\mapsto \frac\sp{(1+\sp x)^2}\mathbbm1_{\R_+}(x)$, along with \cite[2.1.6]{David-Nagaraja}, we get that the density function of $\Tnk nk$ under $\Pninf$ is $s\mapsto \sp (n-k)\binom {n-1}{n-k}  \frac{(\sp s)^{n-k-1}}{(1+\sp s)^n}\mathbbm1_{s\geq0}$. Then \begin{align*}
\E_n^{(\infty)}((\Tnk nk) ^m)&=\sp^{-m}(n-k)\binom {n-1}{n-k} \int_0^\infty \frac{s^{n+m-k-1}}{(1+s)^n}\d s.
\end{align*}
We conclude using Proposition \ref{prop_app_J} in the Appendix.
\end{demopr}

\par\bigskip
\begin{demoth}{\ref{th_loi_CPP_Tor_fixe}}
For any $(t_1,\ldots,t_{n-1})\in(\R_+)^{n-1}$, write
\begin{multline*}
 \P^t(\Hetoile_1\leq t_1,\ \ldots\ ,\Hetoile_{n-1}\leq t_{n-1},\ I_n\leq \mathcal N \leq I_{n+1}\sachant\mathcal N\geq1) \\
  =\sum_{k_0,\ldots,k_n\geq1} \P^t(\Hetoile_1\leq t_1,\ \ldots\ ,\Hetoile_{n-1}\leq t_{n-1},\ I_n\leq \mathcal N \leq I_{n+1},\ I_1=k_0,\ldots,I_n=k_0+\ldots +k_n\sachant\mathcal N\geq1).
\end{multline*}
Now recall from \cite[Th.5.4]{ALContour} that conditional on $\mathcal N\geq 1$, the sequence $(H_i)_{1\leq i\leq \mathcal N-1}$ is distributed as a sequence of i.i.d. random variables satisfying $\P^t(H_i\leq u)=\frac{Nu}{1+Nu}$, stopped at the first one exceeding $t$. Remembering that $\Hetoile_i=\max\{{H_{I_i+1}},\ldots,{H_{I_{i+1}}}\}$, and from the definition of the sequence  $(I_i)_{i\geq1}$,
\begin{align*}
&\P^t(\Hetoile_1\leq t_1,\ \ldots\ ,\Hetoile_{n-1}\leq t_{n-1},\ I_n\leq \mathcal N \leq I_{n+1}\sachant\mathcal N\geq1) \\
 & =\sum_{k_0,\ldots,k_n\geq1} \Bigg(\prod_{i=0}^{n}\frac{\sp }{N}\left(1-\frac{\sp }{N}\right)^{k_i-1} \\
 &\qquad\qquad\qquad\times\P^t(\underset{1\leq i<l_0}\max H_i\leq t,\ \underset{l_0\leq i<l_1}\max H_i\leq t_1,\ldots,\ \underset{l_{n-2}\leq i<l_{n-1}}\max H_i\leq t_{n-1},\ \underset{l_{n-1}\leq i<l_{n}}\max H_i> t) \Bigg) \\
&= \sum_{k_0,\ldots,k_{n}\geq1} \left[\,\prod_{i=0}^{n} \frac{\sp }{N}\left(1-\frac{\sp }{N}\right)^{k_i-1}\right] \left(\frac{Nt}{1+Nt}\right)^{k_0-1} \left[1-\left(\frac{Nt}{1+Nt}\right)^{k_{n}}\right] \prod_{i=1}^{n-1}\left(\frac{N(t_i\wedge t)}{1+N (t_i\wedge t)}\right)^{k_i},
\end{align*}
where for any $0\leq i\leq n$, $l_i:=k_0+\ldots+k_i$.\\

 Now $\forall u\in\R_+$, 

$$ \sum_{k\geq1} \frac \sp  N \left(1-\frac{\sp }{N}\right)^{k-1} \left(\frac{Nu}{1+Nu}\right)^{k}=\frac{\sp u}{1+\sp u},$$

 and
 $$ \sum_{k_{n}\geq1} \frac \sp  N \left(1-\frac{\sp }{N}\right)^{k_{n}-1} \left(1-\left(\frac{Nu}{1+Nu}\right)^{k_{n}}\right) = \frac{1}{1+\sp u}.$$
 Thus we have 
 \begin{multline}\label{eq1}
 \P^t(\Hetoile_1\leq t_1,\ \ldots\ ,\Hetoile_{n-1}\leq t_{n-1},\,I_n\leq \mathcal N \leq I_{n+1} \sachant \mathcal N\geq1) \\
  = \frac{1+Nt}{Nt}\,\frac{1}{1+\sp t}\,\frac{\sp t}{1+\sp t}\prod_{i=1}^{n-1}\frac{\sp (t_i\wedge t)}{1+\sp (t_i\wedge t)}.
 \end{multline}

Finally, by taking $t_i=t$ for all $1\leq i\leq n-1$ in \eqref{eq1}, we get
\begin{equation}
 \label{eq2}
\P^t( I_n\leq \mathcal N \leq I_{n+1} \sachant \mathcal N\geq1) =\frac{1+Nt}{Nt}\frac{1}{1+\sp t}\left(\frac{\sp t}{1+\sp t}\right)^n.
\end{equation}
As a consequence, we have for any $(t_1,\ldots,t_{n-1})\in(\R_+)^{n-1}$,

$$\Pnt(\Hetoile_1\leq t_1,\ \ldots\ ,\Hetoile_{n-1}\leq t_{n-1})=\left(\frac{1+\sp t}{\sp t}\right)^{n-1}\ \prod_{i=1}^{n-1}\frac{\sp(t_i\wedge t)}{1+ \sp(t_i\wedge t)},$$

which leads to the announced result.
\end{demoth}
\par\bigskip

\subsection{Random time of origin}\label{sec_Tor_prior}

We now want to randomize the time of origin. To this aim, we give a (potentially improper) prior distribution to the time of origin in the model defined above. We investigate here priors with density function $g_i:u\mapsto u^{-i}\mathbbm1_{\R_+^*}(u)$, $i\in\Z_+$. The case $i=0$ (resp. $i=1$) is usually referred to as uniform (resp. log-uniform) prior on $(0,\infty)$. \\

For any $0\leq i<n$, recall that $\Pni$ denotes the law of the rescaled tree $\Tree$, with prior $g_i$ on its time of origin, and conditioned on having $n$ sampled individuals at present time :

$$\Pni(\Tree\in\point)=\frac{\int_0^{+\infty} \P^t_n(\Tree\in\point)\P^t(I_n\leq \mathcal N \leq I_{n+1}) g_i(t)\,\d t}{\int_0^{+\infty}\P^{t}(I_n\leq \mathcal N \leq I_{n+1}) g_i(t)\,\d t}.$$

Note that we would have obtained the same distribution $\Pni$ if we had randomized the time of origin before having rescaled time in the process. 

\begin{proposition}\label{prop_loi_de_T_sous_Pni}
For any $0\leq i<n$, the law of $\Tree$ under $\Pni$ is given by
 \upshape$$\Pni(\Tree\in\point)=\int_0^{+\infty} \P^t_n(\Tree\in\point)\, \hni(t)\,\d t,$$\itshape 
where $$\hni:t\mapsto \sp n \binom{n-1}{i}\frac{(\sp t)^{n-i-1}}{(1+\sp t)^{n+1}}\mathbbm1_{\R_+}(t),$$
i.e., the time of origin of $\Tree$ under $\Pni$ is a random variable \upshape$\Tor$ \itshape with posterior distribution characterized by its probability density function $\hni$.
\end{proposition}

\par\bigskip
\begin{demopr}{\ref{prop_loi_de_T_sous_Pni}}
From \eqref{eq2} and from $\P^t(\mathcal N\geq1)=(1+Nt)^{-1}$, we know that for all $t>0$, $\P^t(I_n\leq \mathcal N \leq I_{n+1}) =\frac1{Nt}\frac{(\sp t)^n}{(1+\sp t)^{n+1}}$. Thus,
\begin{align*}
 \int_0^{+\infty}\P^t(I_n\leq \mathcal N \leq I_{n+1}) g_i(t)\,\d t
&= \frac{\sp^i}N \int_0^{+\infty} \frac{(\sp t)^{n-i-1}}{(1+\sp t)^{n+1}}\,\sp\, \d t 
= \frac{\sp^i}N \frac{1}{(i+1)\binom{n}{i+1}} = \frac{\sp^i}{nN}\binom{n-1}i^{-1},
\end{align*}
using Proposition \ref{prop_app_J} in the Appendix. Finally by definition of $\Pni$,
\begin{align*}
\Pni(\Tree\in\point)&= \frac{Nn}{\sp^i}\binom {n-1}i \int_0^{+\infty} \P^t_n(\Tree\in\point)\  \frac \sp N \frac{(\sp t)^{n-1}}{(1+\sp t)^{n+1}}\,\frac{\d t}{t^{i}} \\
&= \int_0^{+\infty} \P^t_n(\Tree\in\point)\ \sp n \binom{n-1}{i}\frac{(\sp t)^{n-i-1}}{(1+\sp t)^{n+1}} \,\d t,
\end{align*} 
which gives the expected result.
\end{demopr}

\par\bigskip
As a corollary, we have that the genealogy of the sample has the law of the genealogy of a birth-death process with fixed size :
\begin{corollaire}\label{coroll_loi_pin_sous_Pni}
 For any $i\in\Z_+$, the rescaled coalescent point process $\pi_n$ is distributed under $\Pni$ as the coalescent point process of a critical birth-death process with parameter $\sp$, with prior $g_i$ on its time of origin, and conditioned on having $n$ extant individuals at present time. 
\end{corollaire}

\begin{remarque}\label{rem_effet_p_sur_Pni}
From the corollary it is easy to see that the sampling parameter $\sp$ only has a scaling effect on time regarding the distribution of $\pi_n$ under $\Pni$. This remains true under $\Pninf$, but not under $\Pnt$ because of the conditioning on the population size at time $t$ (see Remark \ref{rem_effet_ptN}).
\end{remarque}
\par\bigskip
\begin{democo}{\ref{coroll_loi_pin_sous_Pni}}
 The probability for a critical birth-death process with parameter $\sp$ of having $n$ extant individuals at time $t$ is $\frac{(\sp t)^{n-1}}{(1+\sp t)^{n+1}}$ (see \cite[(1)]{Aldous-Popovic}), hence it differs from $\P^t(I_n\leq \mathcal N \leq I_{n+1})$ only by a factor $\sp /N$, and an easy adaptation of the calculations in the proof of Proposition \ref{prop_loi_de_T_sous_Pni} gives the expected result.
\end{democo}

\par\medskip
Finally we study the moments of the divergence times $(\Tnk nk)_{1\leq k\leq n-1}$. The following proposition states a necessary and sufficient condition for the existence of the $m$-th moment of $\Tnk nk$ under $\Pni$. In the case of a uniform prior ($i=0$), we also recall the explicit formula established in \cite[Cor.2.2]{StadlerNew}.

\begin{proposition}\label{prop_moments_Eni(Tnk)}
 For any $0\leq i<n$, $1\leq k\leq n-1$ and $m\geq1$, the $m$-th moment of $\Tnk nk$ under $\Pni$ is finite iff $m\leq k+i$.\\
Besides, for any  $1\leq k\leq n-1$ and $m\leq k$,
$$\E_n^{(0)}((\Tnk nk) ^m)=\frac{\binom{n-k+m-1}m}{\sp^m \binom km}.$$
\end{proposition}

\begin{demo}
From Theorem \ref{th_loi_CPP_Tor_fixe}, we know that under $\Pnt$, the random variables $(\Hetoile_i)_{1\leq i\leq n-1}$ are i.i.d. Hence we obtain from \cite[2.1.6]{David-Nagaraja} that the random variable $\Tnk nk$, defined as the $k$-th order statistic of the sequence $(\Hetoile_i)_{1\leq i\leq n-1}$, has density function 
$$f_{n,k}^t : s\mapsto \sp (n-k)\binom{n-1}{n-k} \frac{(\sp s)^{n-k-1}}{(1+\sp s)^n} \frac{(1+\sp t)^{n-k}}{(\sp t)^{n-1}}(\sp t-\sp s)^{k-1} \mathbbm1_{s\leq t}$$ under $\Pnt$. As a consequence, we have 
$$\Eni((\Tnk nk) ^m)=\int_0^\infty s^m\left(\int_0^\infty f_{n,k}^t(s) \hni(t) \d t\right) \d s,$$
and then 
\begin{align*}
 \Eni((\Tnk nk) ^m)<\infty \ &\Leftrightarrow\ \int_0^\infty \frac{(\sp s)^{n-k-1+m}}{(1+\sp s)^n}\left(\int_{\sp s}^\infty \frac{(\sp t-\sp s)^{k-1}}{(\sp t)^i\,(1+\sp t)^{k+1}} \d t\right) \d s <\infty \\
&\Leftrightarrow\ \int_0^\infty \frac{s^{n-k-1+m}}{(1+s)^n}\left(\int_{s}^\infty \frac{(t-s)^{k-1}}{t^i\,(1+t)^{k+1}} \d t\right) \d s <\infty.
\end{align*}

Let us first characterize the integrability of the function $F:s\mapsto\frac{s^{n-k-1+m}}{(1+s)^n}\left(\int_{s}^\infty \frac{(t-s)^{k-1}}{t^i\,(1+t)^{k+1}} \d t\right)$ in the neighbourhood of $+\infty$. We prove here that $\int_{s}^\infty \frac{(t-s)^{k-1}}{t^i\,(1+t)^{k+1}} \d t\underset{s\to+\infty}\sim cs^{-i-1}$, where $c$ is a (positive) constant w.r.t. $s$. Expanding $(t-s)^{k-1}$, we have
$$\int_{s}^\infty \frac{(t-s)^{k-1}}{t^i\,(1+t)^{k+1}} \d t = \sum_{j=0}^{k-1} (-s)^{k-1-j} \int_s^\infty \frac{\d t}{t^{i-j}(1+t)^{k+1}}.$$
Noting that for any $0\leq j\leq k-1$,
$$\frac{1}{(k+i-j)(1+s)^{k+i-j}}\leq \int_{s}^\infty \frac{\d t}{t^{i-j}(1+t)^{k+1}}\leq \frac{1}{(k+i-j)s^{k+i-j}},$$ 
we obtain 
$$\sum_{j=0}^{k-1} \frac{(-1)^{k-1-j}}{k+i-j}\binom{k-1}j \frac{s^{k+i-j}}{(1+s)^{k+i-j}} \leq s^{i+1}\int_{s}^\infty \frac{(t-s)^{k-1}}{t^i\,(1+t)^{k+1}} \d t \leq \sum_{j=0}^{k-1} \frac{(-1)^{k-1-j}}{k+i-j}\binom{k-1}j,$$
Letting $s\to\infty$ leads to the announced equivalent. As a consequence, $F(s)\underset{+\infty}\sim cs^{m-k-i-2}$, and $F$ is integrable in the neighbourhood of $+\infty$ iff $m-k-i\leq0$.\\

On the other hand, in the case $m-k-i\leq0$, the integrability of $F$ on any compact set of $\R_+$ is clear. Thus $\Eni((\Tnk nk) ^m)$ is finite iff $m-k-i\leq0$. 
%
\end{demo}

\section{Expected frequency spectrum}\label{sec_mut}

\subsection{Mutation setting}

\begin{figure}
 \begin{center}
  \includegraphics[width=.8\linewidth]{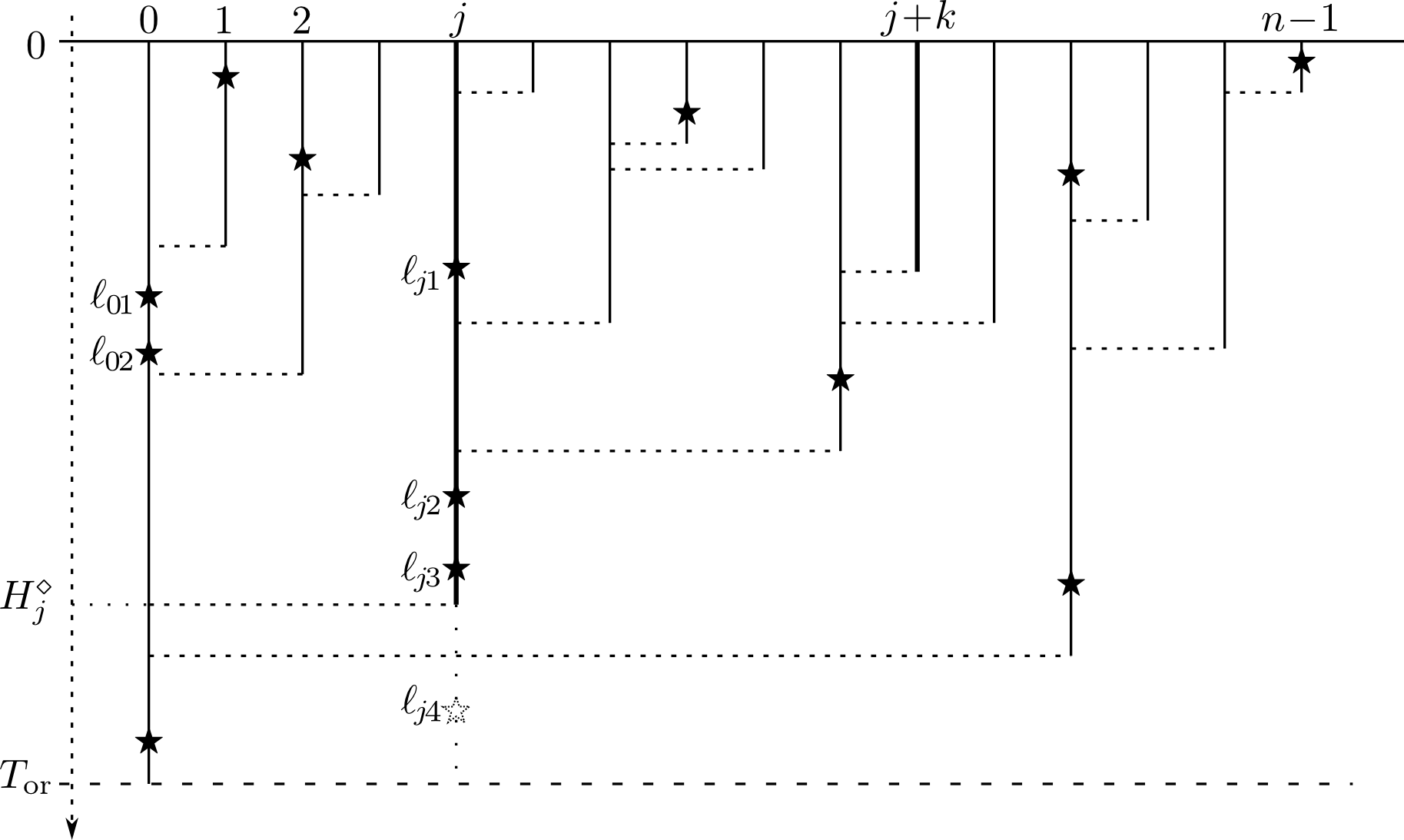}
  \caption{The coalescent point process of a sample of size $n$, with mutations symbolized by stars. Mutations $\ell_{j2}$ and $\ell_{j3}$ are carried by individual $j+k$ while mutation $\ell_{j1}$ is not. Since $\ell_{j4}>\Hetoile_j$, it is not considered as a mutation event. Only mutations $\ell_{01}$, $\ell_{02}$ and $\ell_{j1}$ are carried by two individuals, so that here $\xi_2=3$.}
  \label{fig_mut}
 \end{center}
\end{figure}

Recall from Section \ref{sec_new_intro} that we assume Poissonian mutations at rate $\theta\in\R_+$ on the lineages. We adopt the notation introduced in \cite{ALallelicpartition}, whose framework is very close to ours. Let $(\mathcal P_j)_{j\in\{0,\ldots,n-1\}}$ be independent Poisson measures on $\R_+^*$ with parameter $\theta$. For each $j$ we denote the atom locations of $\mathcal P_j$ by $\ell_{j1}<\ell_{j2}<\ldots$. The branch lengths $(\Hetoile_0,\Hetoile_1,\ldots,\Hetoile_{n-1})$, where we set $\Hetoile_0:=\Tor$, characterize the genealogy of the $n$ individuals (labeled accordingly from 0 to $n-1$) jointly with the foundation time of the population. Then the times $\ell_{jl}$ satisfying $\ell_{jl}<\Hetoile_j$ are interpreted as mutation events, and for all $k\in\{0,\ldots,n-j-1\}$, individual $j+k$ bears mutation $\ell_{jl}$ if 

$$\max\{\Hetoile_{j+1},\ldots,\Hetoile_{j+k}\}<\ell_{jl}<\Hetoile_j,$$

where $\max\varnothing=0$ (see Figure \ref{fig_mut}). The first inequality expresses the fact that a mutation on branch $j$ in the coalescent point process is carried by individual $j+k$ if the time at which it appears is greater than the divergence time of individuals $j$ and $j+k$ (recall that time is running backwards). The second inequality means that all the values $\ell_{jl}$ that are greater than the $j$-th node depth $\Hetoile_j$ are not taken into account.\\

For any $k\in\{1,\ldots,n-1\}$, recall that we denote by $(\xi_k)_{1\leq k\leq n-1}$ the site frequency spectrum of the sample, i.e. $\xi_k$ is the number of mutations carried by $k$ individuals among the $n$ sampled individuals. The sum $S=\sum_{k=1}^{n-1} \xi_k$ is the so-called number of \textit{polymorphic sites}, also known as \textit{single nucleotide polymorphisms} in population genomics. 

\subsection{Results}\label{sec_mut_results}

In this section we give explicit formulae for the expected site frequency spectrum in the case of a fixed time of origin and in the case of a uniform or log-uniform prior on the time of origin. The proofs are based on two different methods, depending on the assumption on $\Tor$, and are expanded in the next section.

\subsubsection{Fixed (finite) time of origin}

The expected site frequency spectrum of the $n$ sampled individuals under $\Pnt$ is given by

\begin{proposition}\label{prop_Ent(xi_k)}
 For any $k\in\{1,\ldots,n-1\}$, $t\in\R_+^*$, defining $\tau:=\sp t$, we have
\begin{multline*}
 \E_n^t(\xi_k)
=\frac\theta\sp \ \Bigg\{\frac{n-3k-1}{k}+\frac{(n-k-1)(k+1)}{k\tau} \\
 \quad+\frac{(1+\tau)^{k-1}}{\tau^{k+1}}\Big[2\tau^2-(n-2k-1)2\tau-(n-k-1)(k+1)\Big]\bigg[\ln(1+\tau)-\sum_{i=1}^{k-1}\frac1i\left(\frac{\tau}{1+\tau}\right)^{i}\bigg] \Bigg\}.
\end{multline*}
\end{proposition}

\subsubsection{Infinite time of origin}

The following two propositions are direct consequences of Proposition \ref{prop_Ent(xi_k)}. However note that Proposition \ref{prop_Eninfty(xi_k)} can be proved independently from the formula provided by Proposition \ref{prop_Ent(xi_k)}, as will be explained in Section \ref{sec_mut_proofs}.

\begin{proposition}\label{prop_Eninfty(xi_k)}
For any $k\in\{1,\ldots,n-1\}$, $\xi_k$ has infinite expectation under $\Pninf$.
\end{proposition}

The infinite expectation of $\xi_k$ under $\Pninf$ leads to consider its renormalization by the expected number of polymorphic sites. The proposition below shows that letting the time of origin go to $+\infty$ flattens the renormalized expected frequency spectrum. A hint for this result is given in Section \ref{sec_mut_proofs}, while we prove it here by letting $t\to\infty$ in Proposition \ref{prop_Ent(xi_k)}.

\begin{proposition}\label{prop_Eninfty(xi_k)_renorm}
For any $k\in\{1,\ldots,n-1\}$,
$$\lim_{t\to\infty}\frac{\E_n^t(\xi_k)}{\E_n^t(S)}=\frac{1}{n-1}.$$
\end{proposition}
\begin{demo}
Fix $k\in\{1,\ldots,n-1\}$. One can easily see from Proposition \ref{prop_Ent(xi_k)} that as $t\to+\infty$, 
$$\E_n^t(\xi_k)\sim 2\theta\ln(t)\ \text{ and }\ \E_n^t(S)\sim 2\theta(n-1)\ln(t),$$
which leads to the result.
\end{demo}

\subsubsection{Random time of origin}

We provide explicit formulae for the expected frequency spectrum for two particular cases of priors : the uniform prior (case $i=0$) and the log-uniform prior (case $i=1$).

\begin{proposition}\label{prop_En0(xi_k)}
 For any $k\in\{1,\ldots,n-1\}$, $\E_n^{(0)}(\xi_k)=n\theta/\sp k$.
\end{proposition}

\begin{proposition} \label{prop_En1(xi_k)}
 For any $k\in\{1,\ldots,n-3\}$, 
$$\E_n^{(1)}(\xi_k)=\frac\theta\sp\,\frac{n(n-1)}{(n-k)(n-k-2)}\Bigg[\frac{n+k-2}{k}-\frac{2(n-1)}{n-k-1}(\mathcal H_{n-1}-\mathcal H_k)\Bigg],$$
where for any $k\in\N$, $\mathcal H_k=\sum_{j=1}^k j^{-1}$.
\end{proposition}

\begin{remarque} 
 The formulae obtained for $\E_n^{(1)}(\xi_{n-2})$ and $\E_n^{(1)}(\xi_{n-3})$, which we chose not to display here, involve non explicit integrals. 
\end{remarque}

\begin{figure}
 \begin{center}
  \includegraphics[scale=.43]{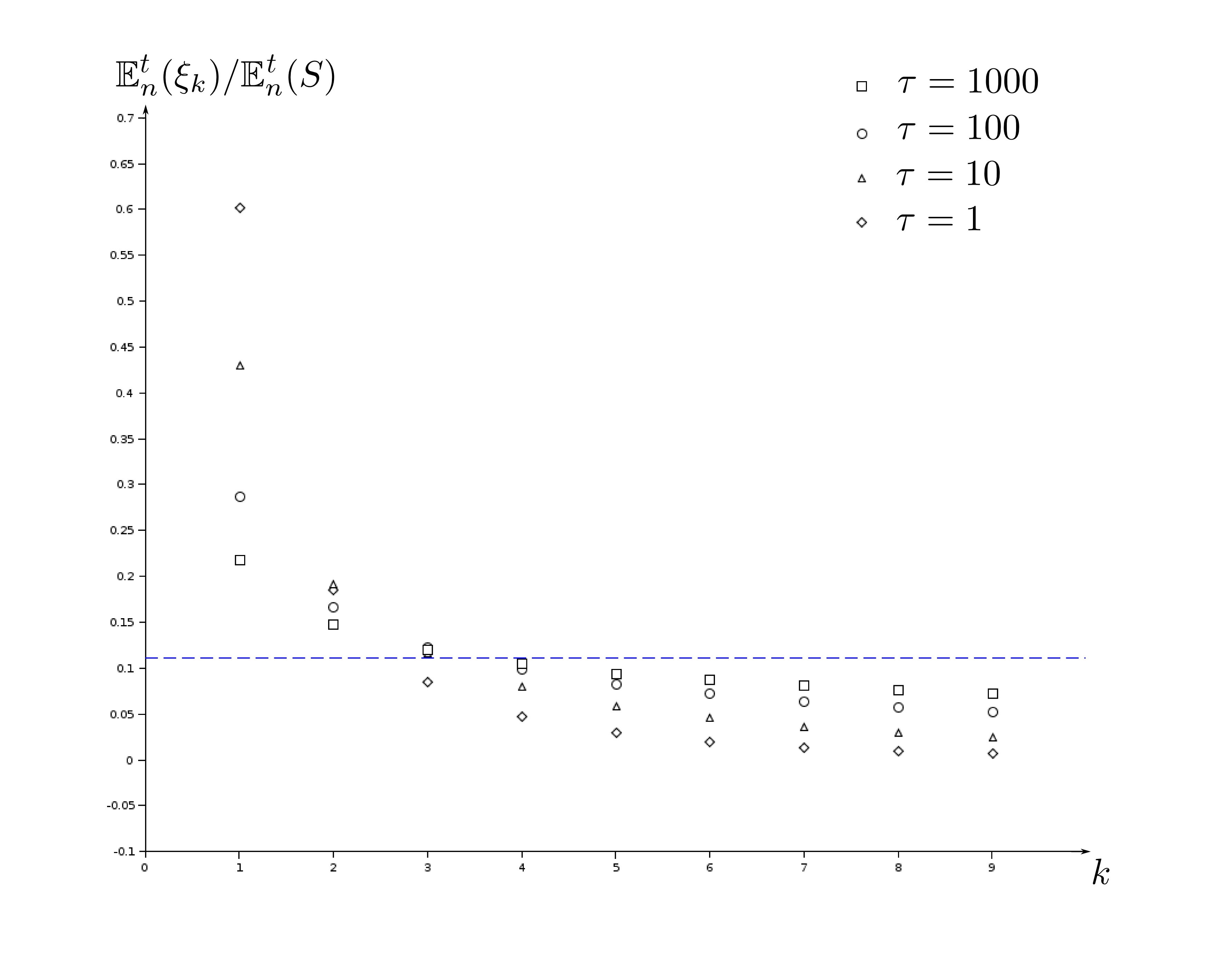}
\caption{The normalized expected site frequency spectrum of a sample of $n=10$ individuals, under $\Pnt$, for $\tau=\sp t\in\{1,10,100,1000\}$. The horizontal dotted line has equation $y=1/(n-1)$}
\label{fig_spt}
 \end{center}
\end{figure}
\begin{figure}
 \begin{center}
  \includegraphics[scale=.43]{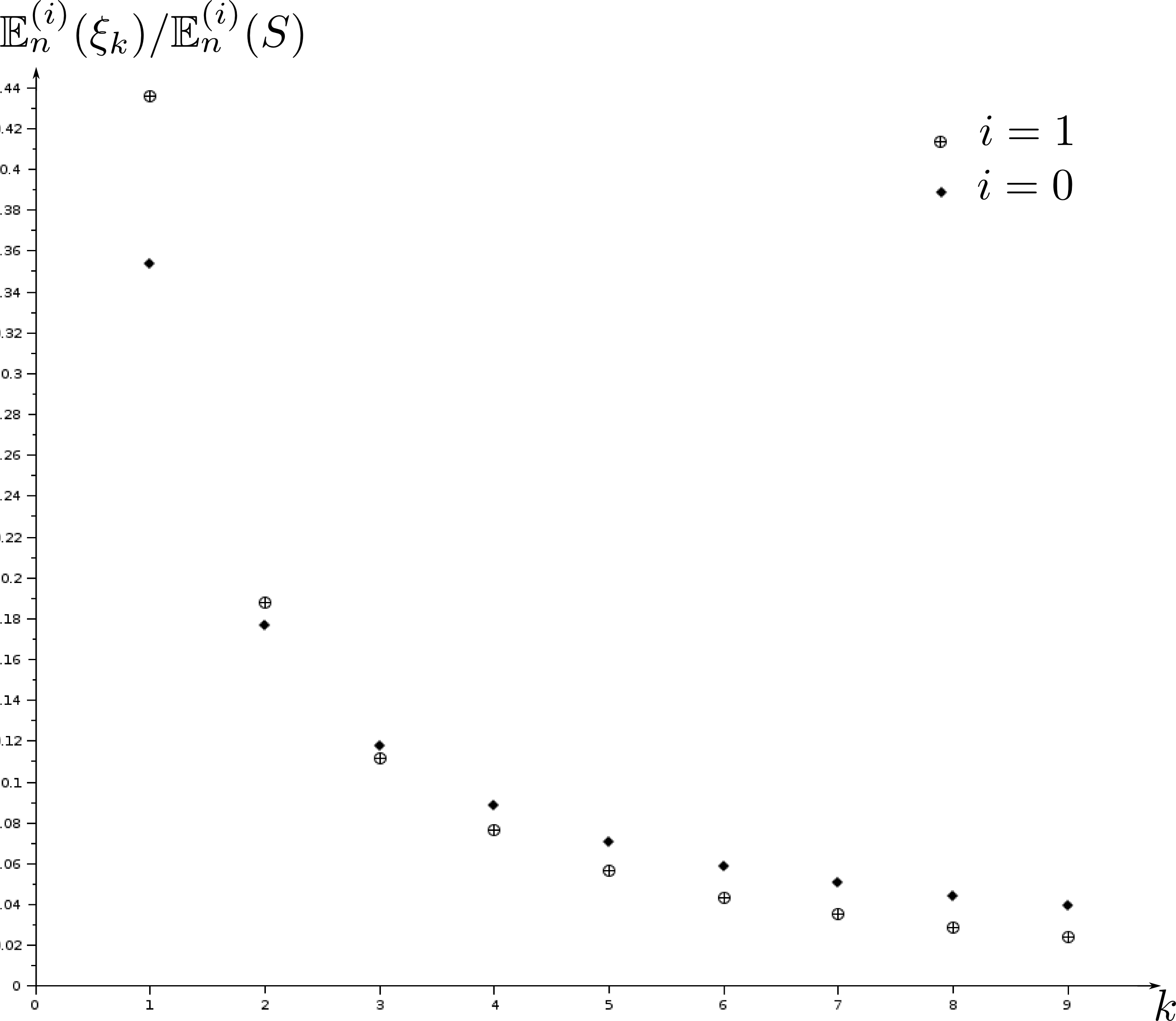}
\caption{The normalized expected site frequency spectrum of a sample of $n=10$ individuals, under $\P_n^{(0)}$ and $\P_n^{(1)}$.}
\label{fig_spp}
 \end{center}
\end{figure}

Graphical representations of the expected frequency spectrum under $\Pnt$, $\P_n^{(0)}$, $\P_n^{(1)}$ are provided in Figures \ref{fig_spt} and \ref{fig_spp}.

\subsection{Proofs}\label{sec_mut_proofs}

Depending on the assumption on $\Tor$, two different methods can be used. The first one relies on an expression of the expected number of mutations carried by $k$ individuals as a function of the expected coalescence times of the tree \cite[pp.103-105]{Wakeley}. The second one decomposes its computation into the sum of the mutations present on lineage $j$, $1\leq j\leq n-k$, carried by $k$ individuals \cite{ALallelicpartition}. Although the second one could be used to prove all the results of Section \ref{sec_mut_results}, the first one provides a very short proof in the cases of an infinite time of origin and of a uniform prior.

\subsubsection{Infinite time of origin and uniform prior}

We base our proof of Propositions \ref{prop_Eninfty(xi_k)} and \ref{prop_En0(xi_k)} on Formula \cite[(4.22)]{Wakeley}, which gives for any $1\leq k\leq n-1$ and any $i\in\Z_+\cup\{\infty\}$

\begin{equation}\label{eq_wakeley}
 \E_n^{(i)}(\xi_k)=\theta\frac 2k \binom{n-1}k^{-1}\  \sum_{j=2}^{n-k+1} \binom j2 \binom{n-j}{k-1} \E_n^{(i)}(\hTnk nj),
\end{equation}
where $\hTnk nj:=\Tnk nj-\Tnk n{j-1}$ denotes the time elapsed between the $(j-1)$-th and the $j$-th coalescence.\\

When the time of origin is set to be infinite a.s., from Proposition \ref{prop_moments_Eninfty(Tnk)} the expected time to the most recent common ancestor is infinite, which entails directly, along with \eqref{eq_wakeley}, that $\E_n^{(\infty)}(\xi_k)$ is infinite for any $k\in\{1,\ldots,n-1\}$. From Equation \eqref{eq_wakeley} we can also give an intuitive explanation of the result of Proposition \ref{prop_Eninfty(xi_k)_renorm}, which establishes that $\lim_{t\to\infty}\E_n^t(\xi_k)/\E_n^t(S)=1/(n-1)$ for any $1\leq k\leq n-1$. Indeed, using \eqref{eq_wakeley} to compute $\Eninf(\xi_k)$, from Proposition \ref{prop_moments_Eninfty(Tnk)} we know that $\Eninf(\hTnk n2)$ is the only infinite contribution to $\Eninf(\xi_k)$. This contribution is thus supported by the first order statistic $\Tnk n1$ of $(\Hetoile_i)_{1\leq i\leq n-1}$ (i.e. the largest divergence time in the coalescent point process). Conditional on $\Tnk n1=\Hetoile_{i_0}$, $\xi_i$ is finite a.s. for any $i\neq n-i_0$. Now under $\Pninf$,
 $(\Hetoile_
i)_{1\leq i\leq n-1}$ is a sequence of i.i.d. random variables, so that the index $i_0$ is uniformly distributed in $\{1,\ldots,n-1\}$. This explains the independence of $\lim_{t\to\infty}\E_n^t(\xi_k)/\E_n^t(S)$ w.r.t. $k$. \\

In the case of a uniform prior on the time of origin, we use a comparison with the very documented Kingman coalescent model. Denote by $\P_{\text{K}}$ the law of the genealogy of a sample of size $n$ under the Kingman coalescent model with mutations at rate $\theta$. First from \cite{StadlerNew} we know that for any $j\in\{2,\ldots,n\}$, the inter-coalescence time $\hTnk nj$ have proportional expectation under $\P_n^{(0)}$ and under the Kingman coalescent model : $\E_n^{(0)}(\hTnk nj)=\frac n{2\sp} \E_{\text{K}}(\hTnk nj)$. Second, from \cite[(4.20)]{Wakeley}, for any $k\in\{1,\ldots,n-1\}$, $\E_{\text{K}}(\xi_k)=\frac{2\theta}k$. As a consequence, using \eqref{eq_wakeley} (which is also valid under $\P_{\text{K}}$) we obtain for any $1\leq k\leq n-1$, $\E_n^{(0)}(\xi_k)=\frac n\sp\frac{\theta}{k}$. This ends the proof of Proposition \ref{prop_En0(xi_k)}.


\subsubsection{Fixed (finite) time of origin and log-uniform prior}

When $\Tor$ is fixed (and finite), or in the case of a non uniform prior on $\Tor$ ($i\in\N$), the equality \eqref{eq_wakeley} does not lead to an explicit expression of the expected frequency spectrum. The formulae stated in Proposition \ref{prop_Ent(xi_k)} (case $\Tor=t\in\R_+^*$) and Proposition \ref{prop_En1(xi_k)} (case of a log-uniform prior on $\Tor$) are obtained using a method developed in \cite{ALallelicpartition} (see proof of Theorem 2.3 for more details).\\

\begin{demopr}{\ref{prop_Ent(xi_k)}}
 Fix $t>0$. Decomposing $\xi_k$ into the sum of the number of mutations on the $j$-th branch carried by exactly $k$ individuals, from \cite{ALallelicpartition} (see proof of Theorem 2.3), we know that
\begin{equation}\label{eq3New}
\E_n^t(\xi_k)=\theta \sum_{j=0}^{n-k} \E_n^t\Big(\big(\min\{\Hetoile_j, \Hetoile_{j+k}\}-\max\{\Hetoile_{j+1},\ldots,\Hetoile_{j+k-1}\}\big)^+\Big), 
\end{equation}
where we have set $\Hetoile_n:=+\infty$. \\

Two particular cases appear, namely $j=0$, where $\min\{\Hetoile_j, \Hetoile_{j+k}\}=\Hetoile_k$ a.s., and $j+k=n$, where $\min\{\Hetoile_j, \Hetoile_{j+k}\}=\Hetoile_j$ a.s. Hence using the i.i.d. property of $(\Hetoile_j)_{1\leq j\leq n-1}$, it follows for any $1\leq j \leq n-k-1$,
\begin{align*}
Q:&=\E_n^t\Big(\big(\min\{\Hetoile_j, \Hetoile_{j+k}\}-\max\{\Hetoile_{j+1},\ldots,\Hetoile_{j+k-1}\}\big)^+\Big)\\
&=\E_n^t\int_0^\infty \mathbbm1_{\max\{\Hetoile_{j+1},\ldots,\Hetoile_{j+k-1}\}<x<\min\{\Hetoile_j, \Hetoile_{j+k}\}} \ \d x\\
&= \int_0^{\infty} \Pnt(\Hetoile_1>x)^2\, \Pnt(\Hetoile_1<x)^{k-1}\ \d x,
\end{align*}

and similarly for $j\in\{0,n-k\}$,
\begin{equation*}
R:=\E_n^t\Big(\big(\min\{\Hetoile_j, \Hetoile_{j+k}\}-\max\{\Hetoile_{j+1},\ldots,\Hetoile_{j+k-1}\}\big)^+\Big)\\
= \int_0^{\infty} \Pnt(\Hetoile_1>x)\, \Pnt(\Hetoile_1<x)^{k-1}\ \d x.
\end{equation*}

From Theorem \ref{th_loi_CPP_Tor_fixe}, we know that $\Pnt(\Hetoile_1<x)=\frac{\sp(x\wedge t)}{1+\sp(x\wedge t)}\frac{1+\sp t}{\sp t}$. This entails, after a change of variables, and recalling that we defined $\tau=\sp t$,

\begin{align*}
 Q&= \frac1\sp\left(\frac{1+\tau}{\tau}\right)^{k-1}\int_0^\tau \left(\frac x{1+x} \right)^{k-1} \left(1-\frac{x(1+\tau)}{\tau(1+x)}\right)^2 \,\d x\\ 
&=\frac1\sp \frac{(1+\tau)^{k-1}}{\tau^{k+1}}\big[\tau^2 I_{k+1,2}(\tau)-2\tau I_{k+1,1}(\tau)+I_{k+1,0}(\tau)\big],
\end{align*}

and
\begin{align*}
R&=\frac1\sp\left(\frac{1+\tau}{\tau}\right)^{k-1} \int_0^\tau \left(\frac x{1+x} \right)^{k-1} \left(1-\frac{x(1+\tau)}{\tau(1+x)}\right) \,\d x  \\
&= \frac1\sp \frac{(1+\tau)^{k-1}}{\tau^{k+1}}\big[\tau^2 I_{k,1}(\tau)-\tau I_{k,0}(\tau)\big],
\end{align*}
where for any $u\in\R_+^*$, $k\in\Z_+$, $l\in\Z$, $I_{k,l}(u):=\int_0^u\frac{x^{k-l}}{(1+x)^k}\d x$. Using Equation \eqref{eq3New}, this leads to

\begin{multline*}
\E_n^t(\xi_k)=\frac\theta\sp \frac{(1+\tau)^{k-1}}{\tau^{k+1}}\Big[(n-k-1)\big(\tau^2 I_{k+1,2}(\tau)-2\tau I_{k+1,1}(\tau)+I_{k+1,0}(\tau)\big) \\ 
+2\big(\tau^2 I_{k,1}(\tau)-\tau I_{k,0}(\tau)\big)\Big] 
\end{multline*}

Finally, using the formulae provided by Proposition \ref{prop_app_I} for $I_{k,l}$, $l\in\{0,1,2\}$, we finally get after some rearrangements
\begin{multline}\label{eq6New}
 \E_n^t(\xi_k) =\frac\theta\sp \,\frac{(1+\tau)^{k-1}}{\tau^{k+1}}\bigg[\ln(1+\tau)\Big(2\tau^2-2(n-2k-1)\tau-(k+1)(n-k-1))\Big) \\
 -2\tau^2+\tau(n-k-1)+\frac{n-k-1}{k}\frac{\tau^{k+2}}{(1+\tau)^k} +(-1)^{k-1}\frac{n-k-1}k\left(1-\frac{1}{(1+\tau)^k}\right)(2\tau+1)\\
\, +\sum_{j=1}^{k-1}\binom{k-1}j \frac{(-1)^j}{j} \left(1-\frac{1}{(1+\tau)^j}\right)\left(2\tau^2+\frac{2\tau k}{j+1}-(n-k-1)\left(2\tau+\frac{k+1}{j+1}\right)\frac{k}{k-j}\right)\bigg].
\end{multline}

To obtain the final form, we decompose the sum in the r.h.s. as follows : 

First define, for $x\in\R$ and $k\in\N$
$$
\begin{array}{ll}
\phii1k(x):=\sum_{j=1}^k\binom kj \frac{x^j}{j}, & \phii2k(x):=\sum_{j=1}^k\binom kj\frac{x^{j+1}}{j(j+1)} \\
\psii1k(x):=\sum_{j=1}^k\binom kj \frac{x^j}{j(k-j)}, \quad & \psii2k(x):=\sum_{j=1}^k\binom kj\frac{x^{j+1}}{j(j+1)(k-j)}.
\end{array}
$$
Then we have
\begin{align*}
 S:=& \sum_{j=1}^{k-1}\binom{k-1}j \frac{(-1)^j}{j} \left(1-\frac{1}{(1+\tau)^j}\right)\left(2\tau^2+\frac{2\tau k}{j+1}-(n-k-1)\left(2\tau+\frac{k+1}{j+1}\right)\frac{k}{k-j}\right) \\
=&\ 2\tau^2\big(\phii1{k-1}(-1)-\phii1{k-1}(-(1+\tau)^{-1})\big) \\
& -2\tau k \big(\phii2{k-1}(-1)-(1+\tau)\phii2{k-1}(-(1+\tau)^{-1})\big) \\
& -(n-k-1)2\tau k \big(\psii1{k-1}(-1)-\psii1{k-1}(-(1+\tau)^{-1})\big)\\
& +(n-k-1)(k+1)k \big(\psii2{k-1}(-1)-(1+\tau)\psii2{k-1}(-(1+\tau)^{-1})\big).
\end{align*}
Let us now reexpress the functions $\phii1k$, $\phii2k$, $\psii1k$ and $\psii2k$. Fix $x\in\R$ and $k\in\N$. The function $\phii1k$ is differentiable at $x$ and we have
$\phii1k'(x)=\sum_{j=1}^k\binom kj x^{j-1} = x^{-1}\big[(1+x)^k-1\big].$ 
This leads by simple integration calculus to $\phii1k(x)=\sum_{j=1}^k \frac{(1+x)^j}{j}-\mathcal H_k$, where $\mathcal H_k=\sum_{j=1}^k j^{-1}$. Then noting that $\phii2k'(x)=\phii1k(x)$, we obtain $\phii2k(x)=\sum_{j=1}^k\frac{(1+x)^{j+1}}{j(j+1)}-x\mathcal H_k+\frac{1}{k+1}-1$. Finally, it is easy to show that $\psii1k(x)=\frac1{k+1}\big(\phii1{k+1}(x)-\frac{x^{k+1}}{k+1}\big)$ and $\psii2k(x)=\frac1{k+1}\big(\phii2{k+1}(x)-\frac{x^{k+2}}{(k+1)(k+2)}\big)$. This yields
\begin{align*}
S=&-2\tau^2 \sum_{i=1}^{k-1}\frac1i\left(\frac\tau{1+\tau}\right)^i\\
  &+2\tau k\left[\frac{k-1} k \tau +\tau\sum_{i=1}^{k-1}\frac1{i(i+1)} \left(\frac\tau{1+\tau}\right)^i \right] \\
  &+2(n-k-1)\tau \left[\sum_{i=1}^{k}\frac1i\left(\frac\tau{1+\tau}\right)^i +\frac{(-1)^{k}}k \left(1-\frac1{(1+\tau)^k}\right)\right] \\
  &+(n-k-1)(k+1)\left[\frac{k}{k+1}\tau -\tau\sum_{i=1}^{k}\frac1{i(i+1)}\left(\frac\tau{1+\tau}\right)^i 
  +\frac{(-1)^{k}}{k(k+1)}\left(1-\frac1{(1+\tau)^k}\right)\right]
\end{align*}
\begin{multline*}
=(n-k-1)k\tau - 2(k-1)\tau^2 +\frac{n-k-1}{k}\frac{\tau^{k+1}}{(1+\tau)^k}+(-1)^{k}\frac{n-k-1}{k} \left(1-\frac1{(1+\tau)^k}\right)(1+2\tau)\\
+(2\tau(n-k-1)-2\tau^2)\sum_{i=1}^{k-1}\frac1i\left(\frac\tau{1+\tau}\right)^i+(2\tau^2k-(n-k-1)(k+1)\tau)\sum_{i=1}^{k-1}\frac1{i(i+1)}\left(\frac\tau{1+\tau}\right)^i\\
= 2\tau^2 -\tau(n-k-1)+\frac{n-k-1}{k}\frac{\tau^{k+1}}{(1+\tau)^k}+(-1)^{k}\frac{n-k-1}{k} \left(1-\frac1{(1+\tau)^k}\right)(2\tau+1)\\
-\frac1k \left(\frac\tau{1+\tau}\right)^{k-1}\left(2k\tau^2-(n-k-1)(k+1)\tau\right)\\
+\left(2(n-2k-1)\tau-2\tau^2+(k+1)(n-k-1)\right)\sum_{i=1}^{k-1}\frac1i\left(\frac\tau{1+\tau}\right)^i,
\end{multline*}
where the last equality was obtained by writing $\frac1{i(i+1)}=\frac1i-\frac1{i+1}$. 
It suffices now to reinject this formula into equation \eqref{eq6New} to obtain the announced result. 
\end{demopr}

\par\bigskip
\begin{demopr}{\ref{prop_En1(xi_k)}}
Reasoning as in the proof of Proposition \ref{prop_Ent(xi_k)}, we express $\E_n^{(1)}(\xi_k)$ as 
\begin{equation}\label{eq4}\E_n^{(1)}(\xi_k)=\theta\, \sum_{j=0}^{n-k} \E_n^{(1)}\Big(\big(\min\{\Hetoile_j, \Hetoile_{j+k}\}-\max\{\Hetoile_{j+1},\ldots,\Hetoile_{j+k-1}\}\big)^+\Big),\end{equation}
with for any $1\leq j \leq n-k-1$,
\begin{align*}
Q:&=\E_n^{(1)}\Big(\big(\min\{\Hetoile_j, \Hetoile_{j+k}\}-\max\{\Hetoile_{j+1},\ldots,\Hetoile_{j+k-1}\}\big)^+\Big)\\
&= \int_0^{\infty} \left( \int_0^\infty  h_n^{(1)}(\tau)\, \Pnt(\Hetoile_1>x)^2\, \Pnt(\Hetoile_1<x)^{k-1}\ \d \tau \right)\, \d x,
\end{align*}

and for $j\in\{0,n-k\}$,
\begin{align*}
R:&=\E_n^{(1)}\Big(\big(\min\{\Hetoile_j, \Hetoile_{j+k}\}-\max\{\Hetoile_{j+1},\ldots,\Hetoile_{j+k-1}\}\big)^+\Big)\\
&=\int_0^{\infty} \left( \int_0^\infty  h_n^{(1)}(\tau)\, \Pnt(\Hetoile_1>x)\, \Pnt(\Hetoile_1<x)^{k-1}\ \d \tau \right)\, \d x.
\end{align*}

From Theorem \ref{th_loi_CPP_Tor_fixe}, we know that $\Pnt(\Hetoile_1<x)=\frac{ \sp (x\wedge t)}{1+ \sp (x\wedge t)}\frac{1+\sp t}{\sp t}$, and from Proposition \ref{prop_loi_de_T_sous_Pni},  for all $t\geq0$, $h_n^{(1)}(t)= \sp n(n-1)\frac{(\sp t)^{n-2}}{(1+\sp t)^{n+1}}$. After a change of variables, this leads to

\begin{align*}
 Q&=\frac1\sp n(n-1)\int_0^\infty\left(\frac{x}{1+x}\right)^{k-1} \left(\int_x^\infty \frac{t^{n-k-1}}{(1+t)^{n-k+2}} \left(1-\frac{x(1+t)}{t(1+x)}\right)^2 \,\d t \right) \d x\\ 
&=\frac1\sp n (n-1)\int_0^\infty \frac{x^{k-1}}{(1+x)^{k+1}}\big[J_{n-k+2,3}(x)-2 x\,J_{n-k+2,4}(x)+x^2J_{n-k+2,5}(x) \big] \d x,\\
 R&=\frac1\sp n(n-1)\int_0^\infty\left(\frac{x}{1+x}\right)^{k-1} \left(\int_x^\infty \frac{t^{n-k-1}}{(1+t)^{n-k+2}} \left(1-\frac{x(1+t)}{t(1+x)}\right) \,\d t \right) \d x\\ 
&=\frac1\sp n (n-1)\int_0^\infty \frac{x^{k-1}}{(1+x)^{k}}\big[J_{n-k+2,3}(x)- x\,J_{n-k+2,4}(x) \big] \d x,
\end{align*}
where for any integers $k\geq l\geq2$ and for any positive real number $x$, $J_{k,l}(x):=\int_x^\infty \frac{u^{k-l}}{(1+u)^k}\d u$.\\

Now using \eqref{eq_J} in Proposition \ref{prop_app_J} to express the integrals $J_{k,l}$ in $R$ and $Q$, and using again Proposition \ref{prop_app_J} to calculate the remaining integrals, we obtain for any $k\geq n-3$,
\begin{multline*}
 Q=\frac1\sp\frac{2n(n-1)}{(n-k)(n-k+1)}\bigg[\sum_{j=0}^{n-k-1} \frac{j+1}{(j+k)(j+k+1)(j+k+2)} \\
\qquad\qquad\quad-\frac{2}{(n-k-1)}\sum_{j=0}^{n-k-2}\frac{(j+1)(j+2)}{(j+k+1)(j+k+2)(j+k+3)}\\
 +\frac{1}{(n-k-1)(n-k-2)}\sum_{j=0}^{n-k-3}\frac{(j+1)(j+2)(j+3)}{(j+k+2)(j+k+3)(j+k+4)}\Bigg],
\end{multline*}
$$\ R=\frac1\sp\frac{n(n-1)}{(n-k)(n-k+1)}\bigg[\sum_{j=0}^{n-k-1} \frac{j+1}{(j+k)(j+k+1)} +\frac{1}{n-k-1}\sum_{j=0}^{n-k-2}\frac{(j+1)(j+2)}{(j+k+1)(j+k+2)}\Bigg].$$
Finally, using partial fraction decompositions to calculate the sums in the expressions of $Q$ and $R$, 
\begin{align*}
 Q&=\frac1\sp\frac{n(n-1)}{(n-k)(n-k+1)}\left[\frac{1}{k}+\frac{6}{n-k-2}-\frac{2(2n+k-1)}{(n-k-1)(n-k-2)}(\mathcal H_{n-1}-\mathcal H_k)\right],\\
 R&=\frac1\sp\frac{n(n-1)}{(n-k-1)(n-k+1)}\left[\frac{n+k-1}{n-k}(\mathcal H_{n-1}-\mathcal H_{k-1})-2\right].
\end{align*}

Reinjecting these expressions into equation \eqref{eq4} leads to
$$ \E_n^{(1)}(\xi_k)=\frac\theta\sp[(n-k-1)Q+2R]=\frac\theta\sp\frac{n(n-1)}{(n-k)(n-k-2)}\left[\frac{n+k-2}{k}-\frac{2(n-1)}{n-k-1}(\mathcal H_{n-1}-\mathcal H_k)\right],$$
for any $1\leq k\leq n-3$, which ends the proof.
\end{demopr}
\par\bigskip

\section{Convergence of genealogies in the large sample asymptotic} \label{sec_cv}

In this section we provide convergence results for the distribution of the suitably rescaled genealogy of a sample of size $n$, as $n\to\infty$. Obtaining such asymptotic results requires an additional assumption on the sampling probability : we assume that the sampling parameter $p$ depends on $n$ in such a way that $p=n/\a$, where $\a\in\R_+^*$. This assumption arises naturally, since it ensures that the expected number of sampled individuals is of order $n$. Besides, according to Remark \ref{rem_effet_p_sur_Pni}, note that the parameter $\a$ will only have a scaling effect on time.\\

In the sequel, the symbol $\egalLoi$ means an equality in law, and for any $n>i\geq0$, $\Loi\point\Pni$ refers to the distribution of a random variable or a process under $\Pni$. Finally, $\cvd$ denotes the convergence in distribution. Recall from \cite[Th.16.16]{Kall} that, if $(\gamma_n)$ is a sequence of random measures on $\R^d$ and $\gamma$ a simple point process on $\R^d$, $\gamma_n\cvd\gamma$ iff $\gamma_n(B)\cvd\gamma(B)$ for any compact set $B$ such that $\gamma(\partial B)=0$, where $\partial B$ denotes the boundary of $B$. \\

\subsection{Results} 
\subsubsection*{Convergence of genealogies}
First define, for any $t>0$, $\pi^t$ (resp $\pi$) as the Poisson point measure on $(0,1)\times(0,\a t)$ (resp. $(0,1) \times\R_+^*$) with intensity $\a\d l\, x^{-2}\d x \mathbbm1_{(l,x)\in(0,1)\times(0,\a t)}$ (resp. $\a\d l\, x^{-2}\d x\mathbbm1_{(l,x)\in(0,1)\times\R_+^*}$). \\

Let $(\rho_i)_{i\geq0}$ be a sequence of i.i.d. exponential random variables with parameter $1/\a$, and define for all $i\geq0$ the inverse-gamma random variable $\e_i:=(\rho_0+\ldots+\rho_i)^{-1}$. Then for $i\in\Z_+$, define the pair $(\pii,\Tori)$, where $\Tori$ is a positive random variable, and $\pii$ is a Cox process $\pii$, as
$$\P(\Tori\in\d t,\ \pii\in\cdot )=\P(\e_i\in\d t)\P(\pi^t\in\cdot).$$
In particular, conditional on $\Tori=t$, $\pii$ has the law of the Poisson point measure $\pi^t$.\\

The first theorem states the convergence in distribution of the random measure  $\pi_n$ under $\Pni$, $i\in\Z_+\cup\{\infty\}$. This result is a generalization of Corollary 2 in \cite{Aldous-Popovic}, which provides convergence in distribution of $\pi_n$ under $\P_n^{(0)}$ towards $\pi^{(0)}$. The proof of this convergence, as well as the proof of the generalization we propose, mainly rely on the convergence of $\pi_n$ under $\P_n^t$ towards $\pi^t$, which is established in Theorem 5 in \cite{Popovic}.

\begin{theoreme}\label{th_cv_pin}
We have the following convergences in distribution as $n\to\infty$ :
\begin{enumerate}[\upshape a)\itshape]
 \item \qquad\qquad\qquad\qquad\qquad\qquad\qquad\quad$\Loi{\pi_n}\Pninf \cvd \pi,$
 \item and for any $i\geq0$,  \qquad\qquad\quad\upshape$\Loi{(\pi_n,\Tor)}\Pni \cvd (\pii,\Tori).$
\end{enumerate}
\end{theoreme}
\par\medskip
As a corollary of this theorem, we state the finite dimensional convergence of the divergence times of $\pi_n$ under $\Pni$, $i\in\Z_+\cup\{\infty\}$. We denote by $(T_k)_{k\geq1}$ (resp. $(\Ti k)_{k\geq1}$) the decreasing reordering of the second coordinates of the atoms of $\pi$ (resp. $\pii$). 

\begin{corollaire}\label{coroll_cv_T1_Tk_et_Ti1_Tik}
 Fix $k\in\N$. We have the following convergences in distribution as $n\to\infty$ :
\begin{enumerate}[\upshape a)\itshape]
 \item  \qquad\qquad\qquad\qquad\qquad\quad\upshape $\Loi{(\Tnk n1,\ldots,\Tnk nk)}\Pninf \cvd (T_1,\ldots,T_k),$\itshape
 \item  and for any $i\geq0$,  \quad\upshape $\ \ \Loi{(\Tor,\Tnk n1,\ldots,\Tnk nk)}\Pni \cvd (\Tori,\Ti1,\ldots,\Ti k).$\itshape
\end{enumerate}
\end{corollaire}
\par\medskip

Besides, the limiting distributions appearing in Corollary \ref{coroll_cv_T1_Tk_et_Ti1_Tik} are specified in the following proposition.

\begin{proposition}\label{prop_loi_de_T1_Tk_et_Tori_T1i_Tki} 
 \quad
\begin{enumerate}[\upshape a)\itshape]
\item For any $k\in\N$, the $k$-tuple \upshape$(T_1,\ldots,T_k)$ \itshape is distributed as \upshape $(\e_0,\ldots,\e_{k-1})$\itshape.
 \item For any $i\in\Z_+$, $k\in\N$, the $k+1$-tuple \upshape$(\Tori,\Ti1,\ldots,\Ti k)$ \itshape is distributed as \upshape $(\e_i,\ldots,\e_{i+k})$\itshape.
\end{enumerate}
\end{proposition}

\par\bigskip
The last theorem describes the links between the different random measures obtained in the limit, in Theorem \ref{th_cv_pin}. Before stating this result, let us clarify some definition. Consider $\mu$ any random measure among $\pi$, $\pii$ ($i\in\Z_+$). Conditional on $\mu=\sum_{t\in A} \delta_{(t,y_t)}$, where $A\subset[0,1]$ is a countable set, denoting by $(u,y_u)$ its largest atom, where we refer to the order w.r.t. the second coordinate, we define the random measure $\sum_{t\in A\setminus\{u\}} \delta_{(t,y_t)}$ as the random measure obtained from $\mu$ by removing its largest atom.\\

Proposition \ref{prop_loi_de_T1_Tk_et_Tori_T1i_Tki} establishes in particular that for any $i\in\Z_+$, the time of origin $\Tori$ is distributed as the $(i+1)$-th largest atom of the random measure $\pi$. The following statement is a direct consequence of this result.
\begin{theoreme}\label{th_pi_plus_gd_atome}\quad
For any $i\in\Z_+$, the measure $\pii$ has the distribution of the random measure obtained from $\pi$ by removing its $i+1$ largest atoms. In particular, for any $i\in\N$, the measure $\pii$ has the distribution of the random measure obtained from $\pi^{(i-1)}$ by removing its largest atom.
\end{theoreme}

As a conclusion, in the limit $n\to\infty$, genealogies with different priors on the time of origin can all be embedded in the same realization of the measure $\pi$ : a realization of the limiting coalescent point process with given prior can be obtained by removing from a realization of $\pi$ a given number of its largest atoms.

\subsubsection*{Convergence of the expected site frequency spectrum}

 Recall that mutations are assumed to occur at rate $\theta$ on the lineages. We deduce the following proposition from the results of Section \ref{sec_mut_results}.

\begin{proposition}
 For any $t\in\R_+^*$ and any $i\in\{0,1\}$, for any $k\in\N$ we have
$$\underset{n\to\infty}\lim \E_n^t(\xi_k)=\a\theta/k \ \ \text{ and }\ \ \ \underset{n\to\infty}\lim \E_n^{(i)}(\xi_k)=\a\theta/k.$$
\end{proposition}

In other words, under $\P_n^t$, $\P_n^{(0)}$ and $\P_n^{(1)}$, the expected site frequency spectrum of the sample converges, as the size of the sample gets large, towards the expected frequency spectrum of the Kingman coalescent \cite[(4.20)]{Wakeley}.

\subsection{Proofs}

To begin with, we state the convergence, as $n\to\infty$, of the posterior distribution of the time of origin $\Tor$ under $\Pni$. This result is essential to obtain other convergence results under $\Pni$, since the posterior density function $\hni$ of $\Tor$ is directly involved in the definition of the law $\Pni$. 

\begin{proposition}\label{prop_cv_Tor_Pni}
 For any $i\in\Z_+$, we have the following convergence in law
\upshape$$\Loi \Tor \Pni \cvd \Tori.$$\itshape
\end{proposition}

\begin{lemme}\label{lemme_loi_de_ei_etc}
For any $i\in\Z_+$, the random variable \upshape$\e_i$ \itshape has density function $\hi:t\mapsto \frac{\a^{i+1}e^{-\a/t}}{i!\; t^{i+2}}\mathbbm1_{t>0}$ (i.e. \upshape$\e_i$\itshape\ follows an inverse-gamma distribution with parameters $(i+1,\a)$).
\end{lemme}

\begin{demo}
 Fix $i\in\Z_+$. The random variable $\e_i$ is the inverse of the sum of $i+1$ independent exponential variables with parameter $\a^{-1}$, i.e. the inverse of a Gamma variable with parameters $(i+1,\a^{-1})$. From the known density function $t\mapsto \frac{\a^{i+1} t^{i}e^{-\a t}}{i!}\mathbbm1_{t>0}$ of a $\Gamma(i+1,\a^{-1})$-variable, we deduce that $\e_i$ has density function $t\mapsto \frac{\a^{i+1}e^{-\a/t}}{i!\;t^{i+2}} \mathbbm1_{t>0}$.
\end{demo}

\begin{demopr}{\ref{prop_cv_Tor_Pni}}
Recall first that by definition, $\Tori$ is distributed as $\e_i$, and as a consequence, has density function $\hi$. From Proposition \ref{prop_loi_de_T_sous_Pni}, recalling that $p=n/\a$, the density function of $\Tor$ under $\Pni$ is given by : for all $t>0$, 
\begin{equation}\label{eq5New}
\hni(t)= \frac{n^2}{\a\,i!}(n-1)\ldots(n-i)\left(\frac{nt/\a}{1+nt/\a}\right)^{n+1}\frac{1}{(nt/\a)^{i+2}}= \frac{(n-1)\ldots(n-i)}{n^i}\,\frac{e^{-(n+1)\ln(1+\frac\a{nt})}}{i!\;\a\;(t/\a)^{i+2}}, 
\end{equation}

and hence for all $t>0$
$$\hni(t)\underset{n\to\infty}\to \frac{\a^{i+1}e^{-\a/ t}}{i!\;t^{i+2}}=\hi(t),$$
and the convergence of the density functions $(\hni)$ towards $\hi$ ensures the convergence in law under $\Pni$ of $\Tor$ towards $\Tori$.
\end{demopr}

\par\bigskip

To prove Theorem \ref{th_cv_pin}, we first recall Theorem 5 of \cite{Popovic}, which can be stated as follows.
\begin{lemme}\label{lemme_AP}
 For any $t>0$, as $n\to\infty$, $\Loi{\pi_n}{\P_n^t} \cvd \pi^t. $
\end{lemme}

\begin{demoth}{\ref{th_cv_pin}}
a) From Proposition \ref{prop_loi_CPP_Tor_infini}, under $\Pninf$ the random measure $\pi_n$ is a simple point process with intensity $\sum_{i=1}^{n-1} \delta_{\{i/n\}}(\d l)\,\frac{n\d x}{\a(1+nx/\a)^2}$. As $n\to\infty$, this intensity measure converges weakly towards $\a\d l\, x^{-2}\d x \mathbbm1_{(l,x)\in(0,1)\times(0,+\infty)}$, which is the intensity measure of the Poisson process $\pi$. From \cite[Th.16.18]{Kall}, this is sufficient to prove the convergence in distribution, under $\Pninf$, of $\pi_n$ towards $\pi$.\\

b) Fix $i\in\Z_+$. For any compact set $A$ of $[0,1]\times\R_+$, $B$ Borel set of $\R_+$ of zero Lebesgue measure boundary, and $k\in\Z_+$, we have
$$\Pni(\pi_n(A)=k,\ \Tor\in B)=\int_B \P_n^t(\pi_n(A)=k)\, \hni(t)\,\d t.$$

In order to apply the dominated convergence theorem, we first remark that for all $t>0$, for any $n>i$, 
\begin{equation}\label{eq7New}
 \hni(t)\leq f^{(i)}(t):=\frac{\a^{i+1}}{i!\;t^{i+2}}, 
\end{equation}
as can easily be seen from \eqref{eq5New}. Besides, studying the variations of $\hni$ yields in particular that $\hni$ is a nonnegative function that increases on $\big(0,\a\frac{n-i-1}{n(i+2)}\big)$. Now there exists $\beta>0$ such that for $n$ large enough,  $\a\frac{n-i-1}{n(i+2)}\geq\beta$. Finally  we have from \eqref{eq7New} that for any $n$ large enough and for all $t>0$, 
$$|\hni(t)|\leq f^{(i)}(\beta)\mathbbm1_{t\leq \beta}+f^{(i)}(t)\mathbbm1_{t>\beta},$$
which is integrable on $\R_+$.

It suffices now to invoke Lemma \ref{lemme_AP} and Proposition \ref{prop_cv_Tor_Pni} to deduce by dominated convergence that 
$$\Pni(\pi_n(A)=k,\ \Tor\in B) \underset{n\to\infty}\longrightarrow \int_B \P(\pi^t(A)=k)\, \hi(t)\,\d t=\P(\pii(A)=k,\ \e_i\in B),$$
and this ends the proof.
\end{demoth}

\par\bigskip
\begin{democo}{\ref{coroll_cv_T1_Tk_et_Ti1_Tik}}
Here we only prove a) since the proof of b) is identical. Fix $k\in\N$ and $A_1,\ldots A_k$, Borel sets of $\R_+^*$ of zero Lebesgue measure boundary, satisfying $\sup A_i=\inf A_{i-1}$ for any $i\in\{2,\ldots,k\}$. We set $B_i=(0,1)\times A_i$ for any $1\leq i\leq k$. Then
\begin{multline*}
\ \Pninf(\Tnk n1\in A_1,\ldots,\Tnk nk\in A_k) =\ \Pninf(\pi_n(B_1)=1,\ldots,\pi_n(B_k)\geq 1)\\
\underset{n\to\infty}\longrightarrow \ \P(\pi(B_1)=1,\ldots,\pi(B_k)\geq 1)=\ \,\P(T_1\in A_1,\ldots,T_k\in A_k),
\end{multline*}
where the convergence follows from Theorem \ref{th_cv_pin}. Furthermore, this result clearly still holds if the sets $(A_i)$ satisfy $\sup A_i\leq\inf A_{i-1}$ instead of $\sup A_i=\inf A_{i-1}$.

To obtain the result in the case where $A_1,\ldots A_k$ are non necessarily pairwise disjoint sets, it suffices to to decompose $\cup_{i=1}^k A_i$ into a partition of disjoint Borel sets and to apply the same reasoning as above. Let us prove this in the simple case $k=2$ : 
\begin{align*}
&  \P(\Tnk n1\in A_1,\Tnk n2\in A_2)\\
&\qquad\qquad=\P(\Tnk n1\in A_1\cap A_2,\Tnk n2\in A_1\cap A_2)+\P(\Tnk n1\in A_1\cap A_2,\Tnk n2\in A_2\backslash A_1) \\
&\qquad\qquad+\P(\Tnk n1\in A_1\backslash  A_2,\Tnk n2\in A_1\cap A_2) +\P(\Tnk n1\in A_1\backslash  A_2,\Tnk n2\in A_2\backslash  A_1)\\
&\qquad\qquad\qquad\qquad =\P( \pi_n(A_1\backslash A_2)=0,\ \pi_n(A_1\cap A_2)\geq 2) \\
&\qquad\qquad\qquad\qquad +\P(\pi_n(A_1\backslash A_2)=0,\ \pi_n(A_1\cap A_2)=1,\ \pi_n(A_2\backslash  A_1)\geq1)\\
&\qquad\qquad\qquad\qquad +\P( \pi_n(A_1\backslash A_2)=1,\ \pi_n(A_1\cap A_2)\geq 1) \\
&\qquad\qquad\qquad\qquad +\P(\pi_n(A_1\backslash A_2)=1,\ \pi_n(A_1\cap A_2)=0,\ \pi_n(A_2\backslash  A_1)\geq1),
\end{align*}
and we conclude as above, using Theorem \ref{th_cv_pin}.
\end{democo}

\par\bigskip


\begin{demopr}{\ref{prop_loi_de_T1_Tk_et_Tori_T1i_Tki}}
  a)\quad We base our reasoning on the fact that a Poisson point measure on $\R_+^*$ with intensity measure $\a x^{-2}\d x$ is the pushforward measure by the continuous mapping $x\mapsto x^{-1}$ of a Poisson process with parameter $\a^{-1}$. Let $\nu$ be such a Poisson process. Then for any $a\in\R_+^*$, recalling that $\rho_0$ is an exponential variable with parameter $\a^{-1}$ and $\e_0=\rho_0^{-1}$ a.s.,
  $$\P(T_1\geq a)=\P(\pi((0,1)\times(a,+\infty))\geq1)=\P(\nu(0,a^{-1})\geq1)=\P(\rho_0\leq a^{-1})=\P(\e_0\geq a),$$
  and hence $T_1\egalLoi \e_0$. A similar reasoning shows that for any $k\geq1$, $A_1,\ldots A_k$, Borel sets of $\R_+^*$ of zero Lebesgue measure boundary, satisfying $\sup A_i<\inf A_{i-1}$ for any $i\in\{2,\ldots,k\}$, 
  $$\P(T_1\in A_1,\ldots,\ T_k\in A_k)=\P(\e_0\in A_1,\ldots,\e_{k-1}\in A_k).$$ 
  As in the previous proof, the case where the sets $(A_i)$ are non pairwise disjoint can be proved with the same reasoning, decomposing $\cup_{i=1}^k A_i$ into a partition of disjoint  sets. We can then conclude that $(T_1,\ldots,T_k)$ is distributed as $(\e_0,\ldots,\e_{k-1})$.\\

  b)\quad In the same way, for any $t\in\R_+^*$, a Poisson point measure on $(0,t)$ with intensity measure $\a x^{-2}\d x\mathbbm1_{(0,t)}(x)$ is the pushforward measure by the mapping $x\mapsto x^{-1}$ of the restriction to $(t^{-1},+\infty)$ of a Poisson process with parameter $\a^{-1}$. Then by definition of $\pii$, for any $a,b\in\R_+^*$,
  \begin{align*}
   \P(\Tori\geq b,\ \Ti1\geq a) &= \int_{a\vee b}^{+\infty} \P(T_1^t\geq a)\,\P(\e_i\in\d t)\\
   &= \int_{a\vee b}^{+\infty} \P(\pi((0,1)\times(a,t))\geq1)\,\P(\e_i\in\d t)\\
   &= \int_0^{a^{-1}\wedge \,b^{-1}} \P(\nu((u,a^{-1}))\geq1)\,\P(\e_i^{-1}\in\d u) \\
   &= \P(\nu((\e_i^{-1},a^{-1}))\geq1,\ \e_i^{-1}\leq b^{-1}).
  \end{align*}
  Now for any $i\geq0$, $\e_i^{-1}$ is distributed as the $(i+1)$-th atom of $\nu$, hence 
 $$   \P(\nu((\e_i^{-1},a^{-1}))\geq1,\ \e_i^{-1}\leq b^{-1})=\P(\e_i^{-1}\leq b^{-1},\ \e_{i+1}^{-1}\leq a^{-1}).$$
  As a conclusion, we have $(\Tori,\Ti1)\egalLoi (\e_i,\e_{i+1})$. With a similar reasoning  we obtain the equality in law, for any $k\geq1$, between $(\Tori,\Ti1,\ldots,\Ti k)$ and $(\e_{i},\ldots,\e_{i+k})$.
\end{demopr}

\par\bigskip

\begin{demoth}{\ref{th_pi_plus_gd_atome}}
 We denote by $\bar\pi^{(i)}$ the random point measure obtained from $\pi$ by removing its $i$ largest atoms. By the restriction property of the Poisson point measures, conditional on $T_i=t$, we have $\bar\pi^{(i)}\egalLoi \pi^t$. 
Recalling from Proposition \ref{prop_loi_de_T1_Tk_et_Tori_T1i_Tki}.(i) that $T_i\egalLoi\e_{i-1}$, for any $a\in\R_+^*$ and $k\geq0$ we have 
\begin{align*}
 \P(\bar\pi^{(i)}((a,+\infty))=k)&=\int_0^{+\infty} \P(\bar\pi^{(i)}((a,+\infty))=k \sachant T_i=t)\, h^{(i-1)}(t)\d t\\
&=\int_a^{+\infty} \P(\pi^t((a,t))=k)\, h^{(i-1)}(t)\d t + \mathbbm1_{k=0} \int_0^a h^{(i-1)}(t)\d t\\
&=\P(\pi^{(i-1)}((a,+\infty))=k),
\end{align*}
where the last equality follows from the definition of the Cox process $\pi^{(i-1)}$. 
\end{demoth}

\appendix
\section{Appendix}\label{appendix}

\begin{proposition}\label{prop_app_I}
For any $k\in\N$, $l\in\Z$ satisfying $k\geq l$, and $x\in\R_+$, we define
\upshape $$I_{k,l}(x):=\int_0^t \frac{t^{k-l}}{(1+t)^{k}} \d t$$ \itshape
Then we have
\begin{enumerate}[\upshape(a)]
 \item for $k\geq0$,\quad \upshape $I_{k,0}(x)=\int_0^x \frac{t^{k}}{(1+t)^{k}} \d t = x-k\ln(1+x)+\sum_{j=1}^{k-1} \frac{(-1)^{j-1}}{j}\binom k{j+1} (1-(1+x)^{-j})$\itshape,

 \item for $k\geq1$,\quad \upshape $I_{k,1}(x)=\int_0^x \frac{t^{k-1}}{(1+t)^{k}} \d t = \ln(1+x)+\sum_{j=1}^{k-1} \frac{(-1)^{j}}{j}\binom {k-1}{j} (1-(1+x)^{-j})$ \itshape,

 \item for $k\geq2$,\quad \upshape $I_{k,2}(x)=\int_0^x \frac{t^{k-2}}{(1+t)^{k}} \d t = \frac1{k-1}\big(\frac{x}{1+x}\big)^{k-1}$ \itshape.
\end{enumerate}
\end{proposition}

\begin{demo}
Using the binomial theorem to expand $\frac{t^k}{(1+t)^k}=\big(1-\frac{1}{1+t}\big)^k$, we get
\begin{align*}
 & I_{k,0}(x)=\sum_{j=0}^k \binom kj (-1)^j \int_0^x (1+t)^{-j} \d t,\ \ \ \text{and} \ \ \ I_{k,1}(x)=\sum_{j=0}^{k-1} \binom {k-1}j (-1)^j \int_0^x (1+t)^{-j-1} \d t,
\end{align*}
which easily leads to (a) and (b).
\end{demo}

\begin{proposition}\label{prop_app_J}
 For any $k\in\N$, $l\in\Z$ satisfying $k\geq l$, and any $x\in\R_+$, define \upshape $$J_{k,l}(x):=\int_x^\infty \frac{t^{k-l}}{(1+t)^k} \d t.$$ 
\itshape
 Then for any $t\in\R_+$, $J_{k,l}(t)<\infty$ if and only if $l\geq 2$. In this case we have
\begin{equation}\label{eq_J} J_{k,l}(x)=\sum_{j=0}^{k-l}\frac{x^j}{(1+x)^{j+l-1}}\frac{(j+1)\ldots(j+l-2)}{(k-1)\ldots(k-l+1)},\end{equation}

and in particular  $J_{k,l}(0)= \frac{(l-2)!}{(k-1)\ldots(k-l+1)}=\left[(l-1)\binom{k-1}{l-1}\right]^{-1}.$

\end{proposition}

\begin{demo}

First for any $l\geq2$ and $x>0$, $J_{l,l}(x)=\int_x^\infty \frac{\d t}{(1+t)^l}=\frac{1}{l-1}(1+x)^{1-l}$. 

For any $k\geq l\geq 2$ and $x\geq0$, an integration by parts gives 
$J_{k,l}(x)=\frac{k}{k-l+1} J_{k+1,l} - \frac{x^{k-l+1}}{(k-l+1)(1+x)^k}.$ 
Then, assuming that $J_{k,l}(x)= \sum_{j=0}^{k-l}\frac{x^j}{(1+x)^{j+l-1}}\frac{(j+1)\ldots(j+l-2)}{(k-1)\ldots(k-l+1)}$, we obtain
\begin{align*}
J_{k+1,l}(x) &= \sum_{j=0}^{k-l}\frac{x^j}{(1+x)^{j+l-1}}\frac{(j+1)\ldots(j+l-2)}{k\ldots(k-l+2)}+\frac{x^{k-l+1}}{k(1+x)^k}\\
 &= \sum_{j=0}^{k-l+1}\frac{x^j}{(1+x)^{j+l-1}}\frac{(j+1)\ldots(j+l-2)}{k\ldots(k-l+2)},
\end{align*}
and \eqref{eq_J} is then proved by induction on $k$.
%
\end{demo}

\selectlanguage{english}
\par\bigskip
\nocite{*}
\bibliographystyle{plain}
\bibliography{biblio_new}
\end{document}